\documentclass[11pt]{article}
\usepackage{graphicx}
\usepackage{parskip}

\usepackage{amsfonts}
\usepackage{natbib}

\usepackage{latexsym}

\newtheorem{theorem}{Theorem}[section]

\newtheorem{lemma}{Lemma}[section]
\newtheorem{corollary}{Corollary}[section]

\newtheorem{example}{Example}[section]

\newcommand{\R}{\mathbb{R}}

\newcommand{\Z}{\mathbb{Z}}

\newcommand{\argmin}{\mathop{\arg \min}\limits}

\newcommand{\PP} {{  \rm I\hskip-0.22em P}}

\begin{document}
\centerline {\bf \Large On the conditions }

\centerline {\bf  \Large used to prove oracle results for the Lasso}

\bigskip
\centerline{Sara van de Geer \& Peter B\"uhlmann}

\centerline{ETH Z\"urich}

\bigskip
\centerline{September, 2009}

\begin{abstract} 
Oracle inequalities and variable selection properties for the Lasso in
linear models have been established under a variety of different
assumptions on the design matrix. We show in this paper how the different
conditions and concepts relate to each other.
The restricted eigenvalue condition \citep{bickel2009sal} or the slightly
weaker compatibility condition \citep{vandeGeer:07a} are sufficient for
oracle results.  We argue that  both these conditions allow for a fairly
general class of design matrices. Hence, 
optimality of the Lasso for prediction and estimation holds for more
general situations than what it appears from coherence
\citep{Bunea:07b,Bunea:07c} or restricted isometry
\citep{candes2005decoding} assumptions. 
\end{abstract}


\section{Introduction}\label{introduction.section}

In this paper we revisit some sufficient conditions for oracle inequalities
for the Lasso in regression and examine their relations. Such oracle results
have been derived, among others, by \cite{Bunea:07c}, \cite{geer08},
\cite{zhanghua08}, \cite{meyu09}, \cite{bickel2009sal}, and for the related
Dantzig selector by \cite{candes2007dss} and \cite{koltch09b}. Furthermore,
variable selection 
properties of the Lasso have been studied by \cite {Meinshausen:06},
\cite{Zhao:06}, \cite{lounici08}, \cite{zhang09} and \cite{wainwright09}. 
Our main aim is to present an overview of the relations (of which some are
known and some are new), and to emphasize that that sufficient conditions
for oracle inequalities hold in fairly general situations.

The Lasso, which we at first only study in a noiseless situation, is
defined as follows.
Let ${\cal X}$ be some measurable space, $Q$ be a probability measure on ${\cal X}$, 
and
$\| \cdot \|$  be the $L_2 (Q)$ norm. Consider a fixed dictionary of functions
$\{ \psi_j \}_{j=1}^p  \subset L_2 (Q)$, and linear functions
$$f_{\beta} (\cdot) := \sum_{j=1}^p \beta_j \psi_j (\cdot ): \ \beta \in \R^p .$$
Consider moreover a fixed target
$$f^0(\cdot) := \sum_{j=1}^p \beta_j^0 \psi_j(\cdot) . $$
We let $S:= \{ j : \ \beta_j^0 \not= 0 \} $ be its active set, and
$s := |S|$ be the {\it sparsity index} of $f^0$.

For some fixed $\lambda > 0 $,  the Lasso for the noiseless problem is
\begin{equation}\label{lasso-noiseless}
\beta^* := \arg \min_{\beta}\biggl \{  \| f_{\beta} - f^0 \|^2 + \lambda \|
\beta \|_1\biggr \} ,
\end{equation}
where $\| \cdot \|_1$ is the $\ell_1$-norm.
We write $f^*: = f_{\beta^*}$ and let $S_*$ be the active set  of the Lasso.

Let us  precise what we mean by an oracle inequality.
With
$\beta $ being a vector in $\R^p$, and
${\cal N} \subset \{ 1 , \ldots , p \}$ an index set, we denote by
$$\beta_{j,{\cal N}} := \beta_j {\rm l} \{ j \in {\cal N}\} , \ j=1 , \ldots , p , $$
the vector with non-zero entries in the set ${\cal N}$
(hence, for example $\beta_S^0 = \beta^0$).

{\bf Definition: Sparsity constant and sparsity oracle inequality.} 
{\it The {\rm sparsity constant $\phi_0$} is the largest value $\phi_0>0$
  such that 
Lasso with $\beta^*$ and $f^*$ satisfies the {\rm $\phi_0$-sparsity oracle
  inequality}  
$$\| f^* - f^0 \|^2 + \lambda \| \beta_{S^c}^* \|_1 \le
{ \lambda^2 s \over \phi_0^2 } . $$}

Restricted eigenvalue conditions  (see \cite{koltch09a, koltch09b} 
 and \cite{bickel2009sal})
have been developed to derive lower bounds
for the sparsity constant. We will present these conditions in the next
section. Irrepresentable conditions (see \cite{Zhao:06}) are tailored for
proving variable selection, 
i.e., showing that $S_* = S$, or, more more modestly,
that the symmetric difference $S_* \triangle S$ is small.

\subsection{Organization of the paper}
We start out with, in Section \ref{definitions.section}, an overview of the conditions we will compare, and some pointers to
the literature. Once the conditions are made explicit, we give in
Subsection \ref{summary} a summary
of the various relations. Figure \ref{fig1} displayed there enables to see these
at a single glance.
We give a proof of each of
the indicated (numbered) implications. Sections \ref{eigenregres.section} -
\ref{irrcomp.section} rigorously deal with all the different cases. The
weakest condition is a compatibility condition. Stronger
conditions can rule out many interesting cases.
We illustrate in Section \ref{verify.section}
that one may check compatibility using approximations.
We give several examples, where the compatibility condition holds.
We also give an example where the compatibility condition
yields a major improvement to the oracle result, as compared to
the restricted eigenvalue condition. 
The noisy case,
studied briefly in Section \ref{noise.section},
poses no additional theoretical difficulties. A
lower bound on the regularization parameter $\lambda$ is required, and
implications become somewhat more technical
because all further results depend on this lower bound. Section
\ref{discussion.section} discusses the results.  

\subsection{Some notation}

For a vector $v$, we invoke the usual notation
$$\| v \|_q =\cases{  ( \sum_j |v_j |^q )^{1/q}, & $1 \le q < \infty$ \cr 
\max_{j} | v_j | , & $q = \infty $ \cr } .$$

The Gram matrix is
$$\Sigma:= \int \psi^T \psi d Q ,$$
so that
$$\| f_{\beta}\|^2 = \beta^T \Sigma \beta . $$
  The entries
of $\Sigma$ are denoted by
$\sigma_{j,k} := ( \psi_j , \psi_k)$, 
with $( \cdot , \cdot)$ being the inner product in $L_2 (Q)$.

To clarify the notions we shall use, 
consider for a moment a partition of the form
$$\Sigma := \pmatrix {\Sigma_{1,1} & \Sigma_{1,2} \cr \Sigma_{2,1} & \Sigma_{2,2} \cr},
$$
where $\Sigma_{1,1}$ is an $N \times N$ matrix,
$\Sigma_{2,1}$ is a $(p-N) \times N$ matrix and
$\Sigma_{1,2}:= \Sigma_{2,1}^T$ is its transpose, and where $\Sigma_{2,2} $ is
a $(p-N) \times (p-N)$ matrix.
Such partitions will be play an important role
in the sections to come. 

More generally,
for a set ${\cal N}\subset \{ 1 , \ldots , p \}$ with size $N$, we introduce the $N \times N$ matrix
$$\Sigma_{1,1} ({\cal N} ) := ( \sigma_{j,k} )_{j,k\in {\cal N} } , $$
the $( p- N) \times N $ matrix
$$\Sigma_{2,1} ({\cal N})= ( \sigma_{j,k} )_{j \notin {\cal N}, k \in  {\cal N} } , $$
and the $(p-N) \times (p-N)$ matrix
$$\Sigma_{2,2} ({\cal N} ) := ( \sigma_{j,k} )_{j,k\notin{\cal N} } . $$
We let $\Lambda_{\rm min}^2(\Sigma_{1,1} ({\cal N})) $ be the smallest eigenvalue of
$\Sigma_{1,1} ({\cal N} )$. Throughout, we assume that, for the fixed
active set $S$,  the smallest eigenvalue
$\Lambda_{\rm min}^2 ( \Sigma_{1,1} (S)) $ is strictly positive, i.e., that
$\Sigma_{1,1} (S)$ is non-singular.

We sometimes identify $\beta_{\cal N}$ with the vector $|{\cal N}
|$-dimensional vector 
$\{ \beta_j \}_{j \in {\cal N}} $, and write e.g., 
$$ \beta_{\cal N}^T \Sigma 
\beta_{\cal N} = \beta_{\cal N}^T  \Sigma_{1,1} (\cal N) \beta_{\cal N} . $$

\section{An overview of definitions}\label{definitions.section}

The definitions we will present are conditions on the Gram matrix $\Sigma$,
namely conditions on quadratic forms $\beta^T \Sigma \beta$, where
$\beta$ is restricted to lie in some subset of $\R^p$. 
We first take the set of restrictions
$${\cal R} (L,S):= \{ \beta : \| \beta_{S^c} \|_1 \le L \| \beta_S \|_1 \not= 0 \}  . $$
The compatibility condition we discuss here is from  \cite{vandeGeer:07a}.
Its name is based on the idea that we require the $\ell_1$-norm and
the $L_2(Q)$-norm to be somehow compatible.

{\bf Definition: Compatibility condition.}
{\it We call
$$\phi_{\rm compatible}^2 ( L,S) :=
\min \biggl \{ {s
\| f_{\beta} \|^2  \over \| \beta_S \|_1^2 } :\ \beta \in {\cal R} (L,S) \biggr \} 
 $$
the {\rm $(L,S)$-restricted $\ell_1$-eigenvalue}.\\
The {\rm $(L,S)$-compatibility condition}
is satisfied if $\phi_{\rm compatible} (L,S) >0$
.}

The bound $\| \beta_S \|_1 \le \sqrt {s} \| \beta_S \|_2$ 
(which holds for any $\beta$) leads to two
successively stronger versions of restricted eigenvalues.
We moreover consider supsets ${\cal N}$ of $S$ with size
at most $N$. Throughout in our definitions, $N \ge s$.
We will only invoke $N=s$ and $N=2s$ (for simplicity).

Define the sets of restrictions
$${\cal R}_{\rm adaptive} (L,S) := \{ \beta: \ \| \beta_{S^c} \|_1 \le \sqrt {s} L \| \beta_S \|_2 \} , $$
and for ${\cal N} \supset S$,
$${\cal R} (L,S,{\cal N} ) := \{ \beta \in {\cal R} (L,S): \ \|\beta_{{\cal N}^c} \|_{\infty} \le \min_{j \in {\cal N} \backslash S} | \beta_j | \} ,$$
and
$${\cal R}_{\rm adaptive} (L,S,{\cal N} ) := \{ \beta \in {\cal R}_{\rm adaptive}  (L,S): \ \|\beta_{{\cal N}^c} \|_{\infty} \le \min_{j \in {\cal N} \backslash S} | \beta_j | \} .$$

If $N=s$, we necessarily have ${\cal N} \backslash S = \emptyset$.
In that case, we let
$\min_{j \in {\cal N} \backslash S} | \beta_j | =0$, i.e.,
${\cal R} (L,S,S) = {\cal R} (L,S) $ (${\cal R}_{\rm adaptive} (L,S,S) = {\cal R}_{\rm adaptive} (L,S) $).

The restricted eigenvalue condition is from  \cite {bickel2009sal}
and \cite{koltch09b}. We
complement it with the {\it adaptive} 
restricted eigenvalue condition. The name of the latter is inspired by the
fact that this strengthened version is useful for the development of theory for
the {\it adaptive} Lasso \citep{zou06} which we do not show in this
paper. 

{\bf Definition: (Adaptive) restricted eigenvalue.}
{\it We call
$$\phi^2 (L, S,N) :=  
\min \biggl \{ {
\| f_{\beta} \|^2  \over \| \beta_{\cal N} \|_2^2 } :\  {\cal N} \supset S ,\ | {\cal N} | \le N , \ 
\beta \in {\cal R} (L,S, {\cal N} )   \biggr \} 
 $$
the {\rm $(L,S,N)$-restricted eigenvalue},
and, similarly,
$$\phi_{ {\rm adaptive}} ^2 (L, S,N) :=
\min \biggl \{ {
\| f_{\beta} \|^2  \over \| \beta_{\cal N} \|_2^2 } :\ {\cal N} \supset S, \ |{\cal N} | \le N, \
\beta \in {\cal R}_{\rm adaptive} (L,S, {\cal N} )
 \biggr \} 
 $$
the {\rm adaptive $(L,S,N)$-restricted eigenvalue}.
The  {\rm (adaptive) $(L,S,N)$-restricted eigenvalue condition} holds if $\phi(L,S,N) >0$ 
($\phi_{ \rm adaptive} (L,S,N) >0$)
.}

We introduce the (adaptive) restricted regression condition
to clarify various connections between different assumptions. 

 {\bf Definition: (Adaptive) restricted regression.}
  {\it The {\rm $(L,S,N)$-restricted regression} is
$$\vartheta(L,S,N) := \max \biggl \{ 
{|( f_{\beta_{\cal N}} , f_{\beta_{{\cal N}^c} }) |  \over 
\| f_{ \beta_{\cal N} } \|^2} :\
\ {\cal N} \supset S, \ |{\cal N} | \le N, \
\beta \in {\cal R} ( L, S , {\cal N} )  \biggr \} .$$
The {\rm adaptive $(L,S,N)$-restricted regression} is

\newpage
$$\vartheta_{\rm adaptive}(L,S,N)  :=  
\ \ \ \ \ \ \ \ \ \ \ \ \ \ \ \ \ \  \ \ \ \ \ \ \ \ \ \ \ \ \ \ \ \ \ \ \ \ \ \ \ \ \ \ \ \ \ \ \ \ \ \ \ \ \ \ \ \ $$
$$ \max \biggl \{ 
{|( f_{\beta_{\cal N}} , f_{\beta_{{\cal N}^c} }) |  \over 
\| f_{ \beta_{\cal N} } \|^2} :\ 
\ {\cal N} \supset S, \ |{\cal N} | \le N, \ \beta \in {\cal R}_{\rm adaptive} ( L, S , {\cal N} )
 \biggr \}.$$
The  {\rm (adaptive) $(L,S,N)$-restricted regression condition} holds if $\vartheta(L,S,N) 
<1$ \break
($\vartheta_{ \rm adaptive} (L,S,N) <1$).}

\vskip .1in
Note that ${( f_{\beta_{\cal N}} , f_{\beta_{{\cal N}^c} }) / 
\| f_{ \beta_{\cal N} } \|^2}$ equals the coefficient when regressing
$f_{\beta_{{\cal N}^c} }$ onto $f_{\beta_{\cal N}}$. 

Of course all these definitions depend on the Gram matrix $\Sigma$.
In Sections \ref{verify.section} and \ref{noise.section}, we make this dependence explicit by adding
the argument $\Sigma$, e.g. the $(\Sigma, L, S)$-compatibility
condition, etc.

When $L=1$, the argument $L$ is omitted,
e.g. $\phi_{\rm compatible} (S) := \phi_{\rm compatible} (1,S)$, and e.g., 
the $S$-compatibility condition is then the condition
$\phi_{\rm compatible} (S) > 0$.
The case $L>1$ is mainly needed to handle the situation with
noise, and $L<1$ is of interest when studying the 
{\it adaptive} Lasso (but we do not develop its theory in this paper). 

We now present some definitions from \cite{candes2005decoding}.

{\bf Definition: Restricted orthogonality constant.}
{\it The quantity
$$\theta (S,N) := 
\sup_{{\cal N} \supset S: \ | {\cal N}| \le N } \ 
\sup_{{\cal M} \subset {\cal N}^c , \ |{\cal M} | \le s } 
\sup_{ \beta} 
\biggl |{  ( f_{\beta_{\cal N}}, f_{\beta_{\cal M}})  \over
\| \beta_{\cal N} \|_2 \| \beta_{\cal M} \|_2 } \biggr | , $$
is called the {\rm $(S,N)$-restricted orthogonality constant}.
We moreover define 
$$\theta_{s,N}:= \max \{ \theta (S,N): \ |S|=s \} . $$}

{\bf Definition: Restricted isometry constant.}
{\it The {\rm $N$-restricted isometry constant} is the smallest value
of $\delta_N $ such that for all ${\cal N} $ with
$| {\cal N} | \le N$, 
$$(1- \delta_N ) \| \beta_{\cal N} \|_2^2 \le \| f_{\beta_{\cal N}} \|^2 \le (1+ \delta_N )
\| \beta_{\cal N} \|_2^2 . $$}

{\bf Definition: Uniform eigenvalue.}
{\it The $(S,N)$-uniform eigenvalue is
$$\Lambda^2(S,N) := 
\inf_{{\cal N} \supset S, \ |{\cal N} | \le N } 
\Lambda_{\rm min}^2 ( \Sigma_{1,1} ({\cal N} )) . $$
}

As mentioned before, we always assume that $\Lambda(S,s) > 0$.

 {\bf Definition: Weak restricted isometry.} 
{\it The {\rm weak $(S,N)$-restricted isometry constant} is 
$$\vartheta_{\rm weak-RIP} (S,N) :={\theta(S,N) \over \Lambda^2 (S,N)}. $$
 The {\rm weak $(L,S,N)$-restricted isometry property} holds if
 $\vartheta_{\rm weak-RIP} (S,N)< 1/L$. }
 
 {\bf Definition: Restricted isometry property.} 
{\it The  {\rm RIP constant} is 
$$\vartheta_{\rm RIP} := { \theta_{s,2s} \over
1- \delta_{s} - \theta_{s,s} }. $$
 The {\rm restricted isometry property}, shortly {\rm RIP}, holds if
 $\vartheta_{\rm RIP}< 1$. }

An irrepresentable condition can be found in \cite{Zhao:06}.
We use a modified version which involves only the design but not the true
coefficient vector $\beta^0$ (whereas its sign vector appears in
\cite{Zhao:06}). The reason is that most other conditions considered in this
paper do not depend on $\beta^0$ as well. Our $(L,S,N)$-irrepresentable
condition with $L=1$ and $N=s$ is only slightly stronger than
the condition in \cite{Zhao:06}. 

{\bf Definition: Irrepresentable condition.}\\
 {\it {\bf Part 1.} We call
$$
\vartheta_{\rm irrepresentable} (S,N) :=  
\min_{{\cal N} \supset S: \ | {\cal N} | \le N } \max_{\| \tau_{\cal N}  \|_{\infty} \le 1} 
\| \Sigma_{2,1} ({\cal N}) \Sigma_{1,1}^{-1}  ({\cal N}) \tau_{\cal N}
\|_{\infty} $$ 
the {\rm $(S,N)$-uniform irrepresentable constant}. 
The {\rm $(L,S,N)$-uniform irrepresentable condition}
is met, if $\vartheta_{\rm irrepresentable} (S,N) < 1/L $. \\
{\bf Part 2.}
We say that the {\rm  $(L,S,N)$-irrepresentable condition}
is met, if for some ${\cal N} \supset S$ with $| {\cal N}  |\le N $, and 
all vectors $\tau_{\cal N}$ satisfying  
$\tau_{\cal N}  \in \{ -1 , 1 \}^{|{\cal N}| }$, we have
$$
 \| \Sigma_{2,1} ({\cal N})\Sigma_{1,1}^{-1} ({\cal N}) \tau_{\cal N} \|_{\infty} < 1/L . 
$$
{\bf Part 3.}
We say that the {\rm weak $(S,N)$-irrepresentable condition}
is met, if for 
all
$\tau_{ S}  \in \{ -1 , 1 \}^{s }$, and 
for some ${\cal N} \supset S$ with $| {\cal N}  |\le N $, and for some
$\tau_{{\cal N} \backslash S } \in \{ -1 , 1\}^{|{\cal N} \backslash S|}$, we have
$$
 \| \Sigma_{2,1} ({\cal N})\Sigma_{1,1}^{-1} ({\cal N}) \tau_{\cal N}
 \|_{\infty} \le 1 .  
$$
}

Finally, we present coherence conditions, which are in the spirit
of \cite{Bunea:07b, Bunea:07c}. \cite{cai09b} derive an oracle result under
a tight coherence condition. 

{\bf Definition: Coherence.} {\it The {\rm  $(L,S)$-mutual coherence condition} holds
if
$$\vartheta_{\rm mutual} (S) :=
 { s \max_{j \notin S} \max_{k \in S} | \sigma_{j,k}| \over
 \Lambda^2 ( S,s) } < 1/L . $$
 The {\rm $(L,S)$-cumulative coherence condition} holds if
 $$\vartheta_{\rm cumulative} (S):=
 {\sqrt s \sqrt { \sum_{k \in S} \biggl ( \sum_{j \notin S} | \sigma_{j,k} | \biggr )^2 } 
 \over \Lambda^2 (S,s) } < 1/L . $$}
 
\subsection{Implications for the Lasso and some first relations}

It is shown in  \cite{vandeGeer:07a} that the compatibility condition implies
oracle inequalities for the Lasso. We re-derive
the result for later reference and also for illustrating that the
compatibility condition is just a condition to make the proof  
go through. We also show (again for later reference) the additional $\ell_2$-result  
if one uses the $(S,N)$-restricted eigenvalue condition.

\begin{lemma} \label{oracle.lemma} (Oracle inequality) We have for the
  Lasso in (\ref{lasso-noiseless}),
$$\| f^* - f^0 \|^2 + \lambda \| \beta_{S^c}^* \|_1 \le 
{\lambda^2 s / \phi_{\rm compatible}^2 (S)}. $$
Moreover, letting ${\cal N}_* \backslash S$ being the set of the $N-s$ largest
coefficients $| \beta_j^* |$, $ j \in S^c$, 
$$\| \beta_{{\cal N}_*}^* - \beta_{{\cal N}_*}^0 \|_2^2 \le \lambda^2 s / 
\phi^4 (S, N) .$$
\end{lemma}

{\bf Proof of Lemma \ref{oracle.lemma}.}
The first assertion follows from the Basic Inequality
$$\| f^* - f^0 \|^2 + \lambda \| \beta^* \|_1 \le \lambda \| \beta^0
\|_1, $$
using the definition of the Lasso in (\ref{lasso-noiseless}),
which implies
$$\| f^* - f^0 \|^2 + \lambda \| \beta_{S^c}^* \|_1 \le \lambda \biggl (
 \| \beta^0 \|_1 - \| \beta_S^* \|_1 \biggr ) $$
 $$ \le \lambda \| \beta_S^* - \beta_S^0 \|_1 \le
 \lambda \sqrt {s} \| f^* - f^0 \| / \phi_{\rm compatible} (S) . $$
Note that the last inequality holds because $\beta^* - \beta^0 \in {\cal
  R}(S)$ which follows by its preceding inequality:
\begin{eqnarray*}
\|\beta^*_{S^c}\|_1 = \|\beta^*_{S^c} - \beta^0_{S^c}\|_1 \le \|\beta^*_S -
\beta^0_S\|_1.
\end{eqnarray*}

The second result follows from
 $$ \| \beta_{{\cal N}_*}^* - \beta_{{\cal N}_*}^0 \|_2^2  \le
 \| f^* - f_0 \|^2 /  \phi^2(S, N), $$
and using $\phi_{\rm compatible}(S) \ge \phi(S,N)$.
\hfill $\sqcup \mkern -12mu \sqcap$

\medskip
An implication of Lemma \ref{oracle.lemma} is an $\ell_1$-norm result:
\begin{eqnarray*}
\|\beta^* - \beta^0\|_1 &= &\|\beta^*_{S^c}\|_1 + \|\beta^*_S - \beta^0_S\|_1\\
&\le& \lambda s/\phi^2_{\rm compatible}(S)  + \lambda \sqrt{s} \|f^* -
f^0\|/\phi_{\rm compatible}(S) \\
  &\le & 2
\lambda s / \phi_{\rm compatible}^2(S),
\end{eqnarray*}
where the last inequality is using the first assertion in Lemma
\ref{oracle.lemma}. We also
note that the second assertion in Lemma \ref{oracle.lemma} has most
statistical importance for the case with $N= s$. We will need the case
$N = 2s$ later in our proofs.

\medskip
\cite {Meinshausen:06} and \cite{Zhao:06} prove that the irrepresentable
condition  
is sufficient and essentially necessary for variable selection, i.e., for
achieving $S_* =S$. We will also present a self-contained proof in Section
\ref{irrepresentable.section} where we will show that the
$(S,s)$-irrepresentable condition is sufficient and the weak
$(S,s)$-irrepresentable condition is essentially necessary for variable
selection.  

\cite {bickel2009sal} prove oracle
inequalities under the restricted eigenvalue condition. They assume 
$$ \min \{ \phi(L,S,s): \ |S| =s \}  >0  $$
(where $L$ can be taken equal to one in the noiseless case).

\begin{figure}[!htb]
\begin{center}
\includegraphics[width=5in,height=1.8in]{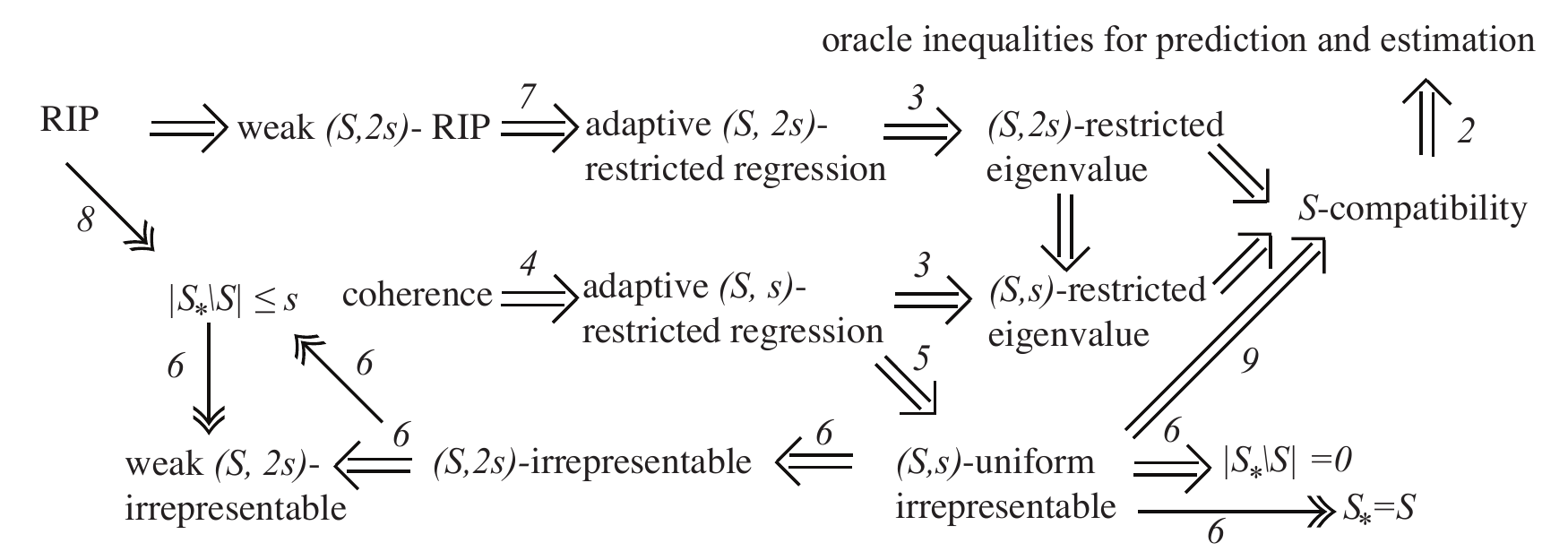}
\end{center}
\caption{A double arrow
($\Rightarrow$) 
indicates a straight implication, whereas the
more fancy arrowheads mean that the relation is under
side-conditions. The numbers indicate the section where the
result is (re)proved. }\label{fig1}
\end{figure}

The restricted isometry property from \cite{candes2005decoding},
abbreviated to RIP, also requires uniformity in $S$. They assume the
RIP 
$$\vartheta_{\rm RIP} < 1 .$$ 
They show that the RIP implies exact reconstruction of $\beta^0$ from $f^0$ by
 linear programming (that is, by minimizing $\| \beta \|_1$ subject to $\|
 f_{\beta} - f^0\| =0$). \cite{cai09a} prove this result assuming $\delta_N +
 \theta_{s,N} < 1$ for $N = 1.25 s$ only; see also \cite{cai09} for an earlier
 result. It is clear that
 $1- \delta_N \le \Lambda^2 (S,N)$,
 i.e., the restricted isometry constants are
 more demanding than uniform eigenvalues.  
 \cite{candes2005decoding} furthermore show that
 $$\vartheta_{\rm weak-RIP} (S,N) \le \vartheta_{\rm RIP} . $$
See also Figure \ref{fig1}. They prove that the RIP is sufficient for
establishing 
oracle inequalities for the Dantzig selector. \cite{koltch09a} and
\cite {bickel2009sal} show  that
$$\phi(L,S,2s) \ge (1-L \vartheta_{\rm weak-RIP} (S, 2s) ) \Lambda (S,2s).$$
Thus, the weak $(S,2s)$-restricted isometry property
implies the $(S,2s)$-restricted eigenvalue condition. See also Figure
\ref{fig1}.  

\cite{Bunea:07a, Bunea:07b, Bunea:07c} show that their coherence conditions
imply oracle results and refinements (see also Section \ref{coherence.section}
for their condition on the diagonal of $\Sigma$). \cite{canpan07} weaken
the coherence conditions by restricting the parameter space for the
regression coefficient $\beta$. 

Finally, it is clear that $ \phi_ {\rm adaptive} ( L,S,N) \le 
\phi(L,S,N) \le \phi_{\rm compatible} (L,S)$, i.e., 

\ \  adaptive restricted eigenvalue condition $\Rightarrow$

\ \  restricted eigenvalue condition
 $\Rightarrow$

\ \ compatibility condition.

See also Figure \ref{fig1}. 

It is easy to see that $\vartheta (L,S,N)$ and
$\vartheta_{\rm adaptive}(L,S,N)$ scale with $L$, i.e., we have
$$\vartheta (L,S,N) = L\vartheta (S,N),\ 
\vartheta_{\rm adaptive} (L,S,N) = L \vartheta_{\rm adaptive} (S,N) . $$
This is not true for the (adaptive) restricted ($\ell_1$-)eigenvalues.
 It indicates that  the (adaptive) restricted regression is not
 well-calibrated for
 proving compatibility or restricted eigenvalue
conditions, i.e, one might pay a large price for taking the route to oracle
results 
via restricted regression conditions.

We end this subsection with the
following lemma, which is based on ideas in \cite{candes2007dss}.
A corollary is the $\ell_2$-bound given in (\ref{l2}), which thus
illustrates that considering supsets ${\cal N}$ of $S$ can be useful.
However, we use the lemma for other purposes as well.

We let for any $\beta$,
 $r_j (\beta):= {\rm rank} ( |\beta_j | )$, $j \in S^c$,
if we put the coefficients in decreasing order.
Let ${\cal N}_0 (\beta) $ be the set of the $s$ largest coefficients in $S^c$:
$${\cal N}_0 (\beta) := \{ j: \ r_j (\beta) \in \{ 1 , \ldots , s \} \}  .$$
Put ${\cal N} (\beta) := {\cal N}_0 (\beta) \cup S$. 
Further, assuming without loss of generality that $p = (K+2)s$
for some integer $K \ge 0$, we let for $k=1 , \ldots , K$, 
$${\cal N}_k (\beta) := \biggl \{ j : \ r_j (\beta)  \in \{ ks +1  , \ldots , (k+1)s \} \biggr \}.$$
We further define
$${\cal N}_* := {\cal N} (\beta^*), \ {\cal N}_k^* := {\cal N}_k (\beta^*), \ k=0, 1 , \ldots , K. $$

\begin{lemma}\label{auxiliary.lemma}
We have for any any $r \ge 1$, and $1/r + 1/q =1$, and any $\beta$,
and for ${\cal N} := {\cal N} (\beta)$, and
${\cal N}_k := {\cal N}_k(\beta)$, $k=0 , 1, \ldots , K$, the bound
$$\| \beta_{{\cal N}^c} \|_r \le \sum_{k=1}^{K} \|\beta_{{\cal N}_k} \|_r \le
\| \beta_{S^c} \|_1 / s^{1/q} . $$
\end{lemma}

\begin{corollary} Combining Lemma \ref{oracle.lemma} with Lemma \ref{auxiliary.lemma} gives
\begin{equation}\label{l2}
\| \beta^* - \beta^0 \|_2^2  \le 2 \lambda^2 s / \phi^4 ( S, 2s) .
\end{equation}
This result is from \cite
{bickel2009sal}. The proof we give is essentially the same as theirs.
\end{corollary}

{\bf Proof of Lemma \ref{auxiliary.lemma}.} 
Clearly,
$$\| \beta_{{\cal N}^c} \|_r = \| \sum_{k=1}^K \beta_{{\cal N}_k} \|_r \le
\sum_{k=1}^{K} \| \beta_{{\cal N}_k } \|_r . $$
We know that for $k=1 , \ldots , K$, 
$$| \beta_j | \le \| \beta_{{\cal N}_{k -1}}  \|_1 / s , \ j \in {\cal N}_k , $$
and hence,
$$\| \beta_{{\cal N}_k } \|_r^r \le s^{-(r-1)}\| \beta_{{\cal N}_{k-1}} \|_1^r  . $$
It follows that
$$\sum_{k=1}^{K} \| \beta_{{\cal N}_k } \|_r \le
\sum_{k=1}^{K} \| \beta_{{\cal N}_{k-1}} \|_1 s^{-(r-1)/r} =
\| \beta_{S^c} \|_1 / s^{1/q}   . $$
\hfill $\sqcup \mkern -12mu \sqcap$

\subsection{Summary of the results}\label{summary}

The following figure summarizes the results. 

Our conclusion is that (perhaps not surprising) the compatibility
condition is the least restrictive, and that many sufficient
conditions for compatibility may be somewhat too harsh
(see also our discussion in Section \ref{discussion.section}).

\section{The restricted regression condition implies the restricted eigenvalue
condition}\label{eigenregres.section}

We start out with an elementary lemma.

\begin{lemma}
Let $f_1$ and $f_2$ by two functions in $L_2(P)$.
Suppose for some $0 <\vartheta < 1 $.
$$- ( f_1 , f_2 )  \le\vartheta \| f_1 \|^2 .$$
Then 
$$ (1-\vartheta) \| f_1\| \le 
\| f_1+f_2 \| . $$
\end{lemma}

{\bf Proof.}
Write the projection of $f_2$ on $f_1$ as
$$f_{2,1}^{\rm P}:= (f_2, f_1 ) / \| f_1\|^2 f_1 . $$
Similarly, let 
$$f=(f_1 + f_2)_1^{\rm P} := (f, f_1 ) / \| f_1\|^2 f_1$$
be the projection of $ f_1 + f_2$ on $f_1$.
Then
$$(f_1 + f_2)_1^{\rm P} = f_1 + f_{2,1}^{\rm P}
= \biggl (1+(f_2,f_1) / \| f_1 \|^2 \biggr ) f_1 ,$$
so that
$$\| (f_1 + f_2)_1^{\rm P}\| =
\biggl |1+(f_2,f_1)/ \| f_1 \|^2 \biggr | \| f_1 \| $$ $$ =
\biggl (1+(f_2,f_1)/ \| f_1 \|^2 \biggr )  \| f_1 \| 
 \ge (1- \vartheta) \| f_1 \| $$
Moreover, by Pythagoras' Theorem
$$ \| f_1+f_2 \|^2 \ge \| (f_1 + f_2)_1^{\rm P} \|^2 .$$
 \hfill $\sqcup \mkern -12mu \sqcap$
 
It is then straightforward to derive the following result. 
\begin{corollary} \label{theta.corollary}
Suppose that  $\vartheta(S,N) < 1/L$.
Then
$$\phi^2(L,S,N) \ge\biggl  (1-
L\vartheta(S,N) \biggr )^2 \Lambda^2 (S,N).$$
A similar result is true for the adaptive versions.
 In other words, the (adaptive) restricted regression condition implies
 the (adaptive) restricted eigenvalue condition.
\end{corollary}

 \section{$S$-coherence conditions imply adaptive $(S,s)$-restricted regression conditions}\label{coherence.section}
 
\cite{Bunea:07a,
Bunea:07b,Bunea:07c} establish oracle results under a condition which we
refer to as the restricted diagonal condition.
They provide coherence conditions for verifying the restricted diagonal
condition. 

{\bf Definition: Restricted diagonal condition.} {\it We say that the {\rm $S$-restricted diagonal condition} holds
if for some constant $\varphi(S)>0$
$$\Sigma - \varphi(S) {\rm diag} (\iota_S) $$
is positive semi-definite.  Here $\iota := ( 1 , \ldots , 1 )^T$
(so $\iota_{j,S} = {\rm l} \{ j \in S\} $).}

We now show that coherence conditions actually imply
restricted regression conditions.
First, we consider some matrix norms in more detail.
Let $1 \le q \le \infty$, and  $r$ be its conjugate, i.e.,
$${1 \over q} + {1 \over r} =1 . $$
Define
$$\| \Sigma_{1,2} ({\cal N}) \|_{2,q}
 := \sup_{\| \beta_{{\cal N}^c}\|_r \le 1} \| \Sigma_{1,2} ({\cal N} ) \beta_{{\cal N}^c} \|_2
 . $$

{\bf Some properties.} {\it The quantity
$\| \Sigma_{1,2} ({\cal N})  \|_{2,2}^2$ is the largest eigenvalue of the matrix
$ \Sigma_{1,2}({\cal N}) \Sigma_{2,1}({\cal N} )$.
We further have for $1 \le q < \infty$, 
$$\| \Sigma_{1,2}({\cal N}) \|_{2,q} \le \left ( \sum_{j\notin {\cal N}} \left ( \sqrt { \sum_{k
\in {\cal N}}
\sigma_{j,k}^2 } \right )^q \right )^{1 / q } 
  ,$$
  and similarly for $q=\infty$,
 $$\| \Sigma_{1,2}({\cal N}) \|_{2,\infty} \le  \max_{j \notin {\cal N} }   \sqrt { \sum_{k
\in {\cal N}}
\sigma_{j,k}^2 }   .$$
   Moreover,
$$
\| \Sigma_{1,2} ({\cal N}) \|_{2,q}\ge \| \Sigma_{1,2} ({\cal N}) \|_{2,\infty},$$
so for replacing
 $\| \Sigma_{1,2} ({\cal N} )\|_{2, \infty} $ by
 $\| \Sigma_{1,2} ({\cal N}) \|_{2, q} $, $q < \infty$, one might have to pay a price. }

\begin{lemma}\label{coherencelemma} 
For all $1 \le q \le \infty$, the following inequality holds:
$$\vartheta_{\rm adaptive}  (S,2s)  \le
    \max_{{\cal N} \supset S,\ |{ \cal N}| =2s }  {  \sqrt s 
 \| \Sigma_{1,2} 
 ({\cal N} ) \|_{2, q } 
 \over  s^{1/q}  \Lambda^2 (S,2s)} 
  . $$
  Moreover,
  $$\vartheta_{\rm adaptive}  (S,s)  \le { \sqrt s \| \Sigma_{1,1} (S) \|_{2,\infty} \over
   \Lambda^2 (S,s)}. $$
  
   \end{lemma}

 \vskip .1in
 {\bf Proof of Lemma \ref{coherencelemma}.} 
 Take $r$ such that $1/q + 1/r =1$.
 Let ${\cal N} \supset S$, with $|{\cal N}| =s$ and let
 $\beta \in {\cal R}_{\rm adaptive}  (S, {\cal N})$. 

 We let $f_{\cal N} := f_{\beta_{\cal N}}$,
 $f_{{\cal N}^c} := f_{\beta_{{\cal N}^c }}$.

We have
$$| ( f_{\cal N} , f_{{\cal N}^c} ) | = | 
 \beta_{\cal N}^T \Sigma_{1,2} ({\cal N} ) \beta_{{\cal N}^c} |  $$
 $$ \le \| \Sigma_{1,2} ({\cal N}) \|_{2,q} \| \beta_{{\cal N}^c }\|_r 
 \| \beta_{\cal N} \|_2 .$$

Applying Lemma \ref{auxiliary.lemma} gives
\begin{equation}\label{rnorm}
\| \beta_{{\cal N}^c} \|_r \le \|  \beta_{S^c} \|_1 / s^{1/q}   \le \sqrt s \| \beta_S \|_2/
s^{1/q} \le \sqrt s \| \beta_{\cal N} \|_2 / s^{1/q} .
\end{equation}
 This yields
$$| ( f_{{\cal N}} , f_{{{\cal N}^c }}  )| \le     
\sqrt s  \| \Sigma_{1,2} (S) \|_{2,q}
\| \beta_{\cal N} \|_2^2 / s^{1/q}$$ $$  \le \sqrt s   \| \Sigma_{1,2} (S) \|_{2,q}  \| f_{\cal N} \|_2^2 / (s^{1/q}
\Lambda^2(S, 2s )).   $$

Similarly,
$$| ( f_{S} , f_{S^c} ) |  
  \le \| \Sigma_{1,2} (S) \|_{2,\infty} \| \beta_{S^c }\|_1 
 \| \beta_S \|_2 $$
$$\le  \sqrt s   \| \Sigma_{1,2} (S) \|_{2,\infty}
 \| \beta_S \|_2^2 \le   \sqrt s   \| \Sigma_{1,2} (S) \|_{2,\infty} / \Lambda^2 (S,s) . $$
\hfill $ \sqcup \mkern -12mu \sqcap$
 
 One of the consequences is in the spirit of the mutual
coherence
condition in \cite{Bunea:07b}.

 \begin{corollary} \label{coherence.corollary1} 
 (Coherence with $q=\infty$) We have
 $$ \vartheta_{\rm adaptive} (S,s) \le 
  {  \sqrt s  \max_{j \notin S}\sqrt   {\sum_{k \in S} \sigma_{j,k}^2 } \over \Lambda^2 (S,s) }
  \le \vartheta_{\rm mutual} (S).$$
   \end{corollary}
   
   With $q=1$ and $N=s$, the coherence lemma is similar to the
cumulative local coherence condition in \cite{Bunea:07c}.
We also consider the case $N=2s$.

\begin{corollary} \label{coherence.corollary3}(Coherence with $q=1$) We have
$$\vartheta_{\rm adaptive} (S,s) \le 
\vartheta_{\rm cumulative} (S), $$and 
$$ \vartheta(S,2s) \le 
\max_{{\cal N} \supset S, \ |{\cal N} | = 2s }  {   \sqrt {  \sum_{k \in {\cal N}}
 \biggl ( \sum_{j \notin {\cal N} } | \sigma_{j,k}| \biggr )^2  }\biggl)  \over
\sqrt s \Lambda^2 (S, 2s)}  . $$
\end{corollary}

The coherence lemma with $q=2$ is a condition about
eigenvalues (recall that $\| \Sigma_{1,2} ({\cal N})\|_{2,2}^2$ equals the largest
eigenvalue of $\Sigma_{1,2} ({\cal N}) \Sigma_{2,1} ({\cal N})$). 
 The bound is then much rougher than the 
one following from the weak $(S, 2s)$-restricted
isometry condition, which we derive in Lemma \ref{RIP.lemma}.
 
 \begin{corollary}\label{coherence.corollary2}(Coherence with $q=2$)
 We have
 $$ \vartheta_{\rm adaptive} (S,2s) \le
  \max_{{\cal N} \supset S,\  |{\cal N} | = 2s } {  \| \Sigma_{1,2} 
 ({\cal N} ) \|_{2, 2 } \over\Lambda^2 (S,2s)
 }   . $$
  \end{corollary}

 \section{ The adaptive $(S,s)$-restricted regression
condition implies the $(S,s)$-uniform irrepresentable condition}
\label{regresirr.section}

\begin{theorem}\label{thetabound.theorem}
We have 
$$\vartheta_{\rm irrepresentable} (S,s)
\le   {\vartheta}_{{\rm adaptive}} (S,s).$$
\end{theorem}

{\bf Proof of Theorem \ref{thetabound.theorem}.}
First observe that
$$  \|\Sigma_{2,1} (S) \Sigma_{1,1}^{-1} (S) \tau_{S}  \|_{\infty}=
\sup_{\| \beta_{S^c} \|_1 \le 1} 
| \beta_{S^c}^T \Sigma_{2,1} (S) \Sigma_{1,1}^{-1} (S) \tau_{S}  | $$
$$ =\sup_{\| \beta_{S^c} \|_1 \le 1} | ( f_{\beta_{S^c}} , f_{b_{S}} ) | , $$
where
$$b_{S} := \Sigma_{1,1}^{-1} (S) \tau_{S} . $$

We note that
$${
\| f_{b_{S}} \|^2 \over \sqrt {s}
\| b_{S} \|_2 }=
{ \| \Sigma_{1,1}^{1/2}  (S) b_{S} \|_2^2 \over 
\| \Sigma_{1,1 } (S)b_{S} \|_2 \| b_{S} \|_2 } 
{\| \Sigma_{1,1} (S) b_{S} \|_2 \over \sqrt {s}} \le 1.$$
(Use Cauchy-Schwarz inequality for bounding the first factor). 
Furthermore, for any constant $c$, 
$$ \sup_{\| \beta_{S^c} \|_1 \le 1}  | ( f_{\beta_{S^c}} , f_{b_{S}} ) |
=  \sup_{\| \beta_{S^c} \|_1 \le c}  | ( f_{\beta_{S^c}} , f_{b_{S}} ) |/c .$$
Take $c= \sqrt {s} \| \beta_{S } \|_2 $ to find
$$  \|\Sigma_{2,1} (S) \Sigma_{1,1}^{-1} (S) \tau_{S}  \|_{\infty}=
 \sup_{\| \beta_{S^c} \|_1 \le  \sqrt {s}  \| b_{S} \|_2  } 
 {  | ( f_{\beta_{S^c}} , f_{b_{S}} ) | \over
  \sqrt {s} \| b_{S} \|_2 }   $$ $$
  \le  \sup_{\| \beta_{S^c} \|_1 \le  \sqrt {s}\| b_{S} \|_2  } 
 {  | ( f_{\beta_{S^c}} , f_{b_{S}} ) | \over
   \|f_{ b_{S}} \|^2 }   .$$

\hfill $\sqcup \mkern -12mu \sqcap$

\section{The $(S,s)$-irrepresentable condition is sufficient and
essentially necessary for variable selection}\label{irrepresentable.section}

An important characterization of the solution $\beta^*$ can be derived from
 the {\it Karush-Kuhn-Tucker} ({\it KKT}) conditions which in our context
 involves subdifferential calculus: see \cite{bertsimas1997introduction}.

{\bf The KKT conditions.} {\it We have 
$$2 \Sigma (\beta^* - \beta^0) = - \lambda \tau^*. $$
Here $\| \tau^* \|_{\infty} \le 1$, and moreover
$$\tau_j^*{\rm l} \{ \beta_j^* \not= 0 \} = {\rm sign} (\beta_j^*)   , \ j = 1 , \ldots , p. $$}

For ${\cal N} \supset S$, we write the projection of a function $f$ on the space spanned by
$\{ \psi_j \}_{j\in {\cal N}}$ as $f^{P_{\cal N}}$, and the anti-projection as
$f^{A_{\cal N}}:= f - f^{P_{\cal N}} $.  Hence, we note that
$$f_{\beta}^{P_{\cal N}} = ( f_{\beta_{\cal N}} + f_{\beta_{{\cal N}^c}} )^{P_{\cal N}} =
f_{\beta_{\cal N}} + (f_{\beta_{{\cal N}^c}})^{P_{\cal N}} , $$
and thus
$$f_{\beta}^{A_{\cal N}} = ( f_{\beta_{{\cal N}^c}} )^{A_{\cal N}} . $$
Moreover
$$\| (f_{\beta_{{\cal N}^c}})^{A_{\cal N}} \|^2 =
\beta_{{\cal N}^c}^T\Sigma_{2,2} ({\cal N})  \beta_{{\cal N}^c}-
  \beta_{{\cal N}^c}^T\Sigma_{2,1} ({\cal N}) \Sigma_{1,1}^{-1} ({\cal N}) \Sigma_{1,2} ({\cal N})
\beta_{{\cal N}^c}. $$

\begin{lemma} \label{start.lemma} Suppose $\Sigma_{1,1}^{-1} ({\cal N})$ exists.
We have
$$2 \| (f_{\beta_{{\cal N}^c}^*})^{A_{\cal N}}\|^2=
\lambda( \beta_{{\cal N}^c}^* )^T
 \Sigma_{2,1}({\cal N} )\Sigma_{1,1}^{-1} ({\cal N}) \tau_{\cal N}^* - \lambda \| \beta_{{\cal N}^c}^* \|_1 . $$
\end{lemma}

{\bf Proof of Lemma \ref{start.lemma}.}
By the KKT conditions, we must have
$$2 \Sigma_{1,1} ({\cal N}) (\beta_{\cal N}^* - \beta_{\cal N}^0 ) +  2\Sigma_{1,2}({\cal N})
\beta_{{\cal N}^c}^* = - \lambda \tau_{\cal N}^* ,$$
$$2 \Sigma_{2,1} ({\cal N}) (\beta_{\cal N}^* - \beta_{\cal N}^0 ) +
2 \Sigma_{2,2} ({\cal N})  \beta_{{\cal N}^c}^* = -\lambda  \tau_{{\cal N}^c}^*.$$
It follows that
 $$2(\beta_{\cal N}^* - \beta_{\cal N}^0 ) + 2\Sigma_{1,1}^{-1} ({\cal N}) \Sigma_{1,2} ({\cal N})
\beta_{{\cal N}^c}^* = - \lambda \Sigma_{1,1}^{-1} ({\cal N}) \tau_{\cal N}^* ,$$
$$2\Sigma_{2,1} ({\cal N}) (\beta_{\cal N}^* - \beta_{\cal N}^0 ) +2
\Sigma_{2,2} ({\cal N})  \beta_{{\cal N}^c}^* = -\lambda \tau_{{\cal N}^c}^*  $$
(leaving the second equality untouched).
Hence, multiplying the first equality  \break by  $-(\beta_{{\cal N}^c}^*)^T\Sigma_{2,1} ({\cal N})$, and
the second by $ (\beta_{{\cal N}^c}^* )^T $,
 $$- 2(\beta_{{\cal N}^c}^* )^T \Sigma_{2,1} ({\cal N}) (\beta_{\cal N}^* - \beta_{\cal N}^0 ) 
 - 2(\beta_{{\cal N}^c}^* )^T\Sigma_{2,1} ({\cal N}) \Sigma_{1,1}^{-1} ({\cal N}) \Sigma_{1,2} ({\cal N})
\beta_{{\cal N}^c}^*  $$ $$=  \lambda (\beta_{{\cal N}^c}^* )^T
 \Sigma_{2,1}({\cal N} )\Sigma_{1,1}^{-1} ({\cal N}) \tau_{\cal N}^* ,$$
$$2(\beta_{{\cal N}^c}^* )^T\Sigma_{2,1} ({\cal N}) (\beta_{\cal N}^* - \beta_{\cal N}^0 ) 
+ 2(\beta_{{\cal N}^c}^* )^T\Sigma_{2,2} ({\cal N})  \beta_{{\cal N}^c}^* = -\lambda \| \beta_{{\cal N}^c}^* \|_1 ,$$
where we invoked that $\beta_j^* \tau_j^* = |\beta_j^* |$.
Adding up the two equalities gives
$$2(\beta_{{\cal N}^c}^* )^T\Sigma_{2,2} ({\cal N})  \beta_{{\cal N}^c}^*-
 2 (\beta_{{\cal N}^c}^* )^T\Sigma_{2,1} ({\cal N}) \Sigma_{1,1}^{-1} ({\cal N}) \Sigma_{1,2} ({\cal N})
\beta_{{\cal N}^c}^* $$
$$ = \lambda (\beta_{{\cal N}^c}^* )^T
 \Sigma_{2,1} ({\cal N})\Sigma_{1,1}^{-1} ({\cal N}) \tau_{\cal N}^* - \lambda \| \beta_{{\cal N}^c}^* \|_1 . $$
 \hfill $\sqcup \mkern -12mu \sqcap$

We now connect the irrepresentable condition to variable
selection. 
Define
$$|\beta^0 |_{\rm min} := \min \{ |\beta_j^0 | : \ j \in S \} . $$

\begin{lemma} \label{select.lemma} $ $\\
{\bf Part 1.} Suppose the $(S,N)$-uniform irrepresentable condition holds.
Then $|S_* \backslash S| \le N-s$.\\
{\bf Part 2.}
Suppose the $(S,N)$-irrepresentable condition holds
and 
$$|\beta_{\rm min}^0 | > \lambda s / \phi_{\rm compatible}^2 (S). $$
Then $S_* \supset S $ and $|S_*| \le N$.\\
{\bf Part 3.}
Conversely, suppose that  $S_* \supset S$ and $|S_* | \le N$, and
$\Lambda(S,N) >0$. Then
$$\| \Sigma_{2,1} (S_*) \Sigma_{1,1}^{-1} (S_*) \tau_{S_*}^* \|_{\infty} \le 1 .$$
If moreover 
$$|\beta^0 |_{\rm min} > \lambda \sqrt {s} / (2 \Lambda (S,N)) , $$
then $\tau_{S_*}^* = \tau_{S_*}^0$, where $\tau_{S_*}^0 := {\rm sign}
(\beta_{S_*}^0)$. 
\end{lemma}

A special case is $N=s$.
In Part 1, we then obtain that $S_* \subset S$, i.e., no false
positive selections. Moreover, Part 2 then proves $S_*=S$ and
Part 3 assumes $S_*=S$. 

{\bf Proof of Lemma \ref{select.lemma}.} 

{\bf Part 1.}
Let ${\cal N} \supset S $ be a set of size  at most $N$, such that
$$
\sup_{\| \tau_S \|_{\infty} \le 1 } \| \Sigma_{2,1} ({\cal N}) \Sigma_{1,1}^{-1} ({\cal N} ) \tau_{\cal N} \|_{\infty} < 1 . $$
By Lemma \ref{start.lemma}, we now have that if $\| \beta_{{\cal N}^c}^* \|_1 >0$
$$2 \| (f^* )^{A_{\cal N}}  \|^2 =
\lambda (\beta_{{\cal N}^c}^* )^T\Sigma_{2,1} ({\cal N}) \Sigma_{1,1}^{-1} ({\cal N} ) \tau_{\cal N}^* - \lambda \| \beta_{{\cal N}^c}^* \|_1  < 0 , $$
which is a contradiction. Hence $\| \beta_{{\cal N}^c}^* \|_1 =0 $, i.e.,
$S_* \subset {\cal N} $.

{\bf Part 2.} 
By Lemma \ref{oracle.lemma},
$$\| \beta_S^* - \beta_S^0 \|_1 \le \sqrt s \| f^*- f^0 \| / \phi_{\rm compatible} (S) \le
\lambda s / \phi_{\rm compatible}^2 ( S) . $$
The condition
$|\beta_{\rm min}^0 | > \lambda {s} / \phi_{\rm compatible}^2 (S)$
thus implies that $S_* \supset S$, and hence that
$\tau_S^* \in \{ -1 , 1 \}^s $. We also know that $\tau_{S_* }^* \in \{ -1 , 1 \} $.
Hence for any ${\cal N}$ satisfying $S \subset {\cal N} \subset S_*$, 
also $\tau_{\cal N} \in \{-1 , 1 \}^{|{\cal N} |} $. Thus, by the $(S,N)$-irrepresentable
condition, there exists such an ${\cal N}$, say $\tilde {\cal N}$, with
$$\| \Sigma_{2,1} (\tilde  {\cal N} ) \Sigma_{1,1}^{-1} (\tilde {\cal N} ) \tau_{\tilde {\cal N}}^* \|_{\infty} < 1 . $$
As in Part 1, we then must have that $\| \beta_{\tilde {\cal N}^c}^* \|_1=0 $.

{\bf Part 3.}
Because $\Lambda(S,N)>0$, and $|S_* | \le N$, we know that $\Sigma_{1,1}^{-1} (S_*) $
exists. Because $S_* \supset S$, we have $\beta_{S_*^c}^* =\beta_{S_*^c}^0 = 0$, 
so the KKT conditions take  the form
$$2 \Sigma_{1,1} (S_*) ( \beta_{S_*}^* - \beta_{S_*}^0 ) =-\lambda \tau_{S_*}^* , $$
and
$$2 \Sigma_{2,1} (S_*) (\beta_{S_*}^* - \beta_{S_*}^0 ) =
-\lambda \tau_{S_*^c}^* . $$
Hence
$$\beta_{S_*}^* - \beta_{S_*}^0 = \lambda \Sigma_{1,1}^{-1} (S_*)\tau_{S_*}^* / 2 , $$
and, inserting this in the second KKT equality,
$$\Sigma_{2,1} (S_*) \Sigma_{1,1}^{-1} (S_*) \tau_{S_*}^* = \tau_{S_*^c}^* . $$
But then
$$\| \Sigma_{2,1} (S_*) \Sigma_{1,1}^{-1} (S_*) \tau_{S_*}^* \|_{\infty}=
\|  \tau_{S_*^c}^*\|_{\infty} \le 1. $$
The first KKT equality moreover implies
$$
\| \beta_{S_*}^* - \beta_{S_*}^0 \|_2 \le \lambda \sqrt {N} / (2\Lambda^2 (S,N)) . $$
So when $|\beta^0 |_{\rm min} > \lambda \sqrt {N} / (2\Lambda^2 (S,N))$, we have
$\tau_{S_*}^* = \tau_{S_*}^0 $.

\hfill $\sqcup \mkern -12mu \sqcap$

\section{The weak $(S,2s)$-restricted isometry property implies 
the $(S,2s)$-restricted regression condition}\label{RIPregres.section}

\begin{lemma} \label{RIP.lemma}We have
$$\vartheta_{{\rm adaptive}} (S,2s) \le \vartheta_{\rm weak-RIP}(S, 2s) . $$
\end{lemma}

{\bf Proof of Lemma \ref{RIP.lemma}.} 
Let $\beta$ be an arbitrary vector.
satisfying $\|\beta_{S^c} \|_1 \le \sqrt s \| \beta_S \|_2 $.
From Lemma \ref{auxiliary.lemma}, 
$$\sum_{k=1}^{K} \| \beta_{{\cal N}_k } \|_2 \le
\| \beta_{S^c} \|_1 / \sqrt s \le \| \beta_{S} \|_2   . $$

Hence, using the definition of the restricted orthogonality
constant $\theta(S,2s)$,  and of the $(S, 2s)$-uniform eigenvalue
$\Lambda^2 (S, 2s)$, 
$$| ( f_{\beta_{\cal N}} , f_{\beta_{{\cal N}^c }} )| \le
\theta(S,2s)  \sum_{k=1}^{K}
\| \beta_{\cal N} \|_2 \| \beta_{{\cal N}_k} \|_2 \le
\theta(S,2s) \| \beta_{\cal N} \|_2 
 \| \beta_S \|_2 $$ $$
\le  \theta(S,2s)   \| f_{\beta_{\cal N}} \|_2^2 /
\Lambda^2 ( S,2s) , $$
or
$${ | ( f_{\beta_{\cal N}} , f_{\beta_{{\cal N}^c}} )| \over
\| f_{\beta_{\cal N}} \|^2}  \le 
\theta (S, 2s)  /
\Lambda^2 ( S, 2s) = \vartheta_{\rm weak-RIP} (S, 2s). $$

\hfill $\sqcup \mkern -12mu \sqcap$

\begin{corollary} \label{RIP.corollary}
Together with Corollary \ref{theta.corollary}, we can now conclude that when
$\vartheta_{\rm weak-RIP}(S, 2s)<1/L$, one has
$$\phi(L,S,2s) \ge (1-L \vartheta_{\rm weak-RIP} )^2 \Lambda^2 (S,2s) . $$
This result is from \cite{koltch09a} and \cite {bickel2009sal}.
\end{corollary}

\section{The restricted isometry property
with small constants implies the weak
$(S,2s)$-irrepresentable condition}\label{RIPirr.section}

We start with two preparatory lemmas.
Recall that 
$$\vartheta_{\rm weak-RIP} (S,s) = \theta(S,s) / \Lambda^2 (S,s)  . $$ 

  \begin{lemma} \label{antiprojection.lemma} Suppose that
 $$\vartheta_{\rm weak-RIP} (S,s)< 1 . $$ 
 Then
$$ 2 \| (f_{\beta_{S^c}^*})^{A_S}\|^2\le \vartheta_{\rm weak-RIP} (S,s)
 \biggl ( \lambda  \sqrt {s}  \| \beta_{{\cal N}_0^*}^* \|_2 \biggr ) , $$
where $A_S$ denotes the anti-projection defined in Section
\ref{irrepresentable.section}. 
 \end{lemma}
 
 {\bf Proof of Lemma \ref{antiprojection.lemma}.}
 Define
 $$b_S:= \Sigma_{1,1} (S)^{-1} \tau_S^* . $$
 Then
 $$\| b_S \| \le \| \tau_S^* \|_2 / \Lambda^2(S,s) \le \sqrt{s} / \Lambda^2(S,s) . $$
 Moreover,
 $$
| ( \beta_{S^c}^* )^T
 \Sigma_{2,1} \Sigma_{1,1}^{-1} (S) \tau_S^* |=
 |( f_{\beta_{S^c}^*} , f_{b_S} ) | \le
 \sum_{k=0}^{K-1} | (f_{\beta_{{\cal N}_k^*}^*} , f_{b_S} ) | $$
 $$ \le \theta (S,s) \sum_{k=0}^{K} \| \beta_{{\cal N}_k^*}^* \|_2 \| b_S \|_2 
  \le \theta (S,s) \| b_S \|_2 \left ( 
 \| \beta_{{\cal N}_0^*}^* \|_2 + \sum_{k=1}^{K}\| \beta_{{\cal N}_k^*}^* \|_2 \right ) 
 $$ $$ \le \theta (S,s) \| b_S \|_2 \left ( 
 \| \beta_{{\cal N}_0^*}^* \|_2 + \| \beta_{S^c}^* \|_1 / \sqrt {s} \right ) 
  \le { \theta (S,s) \over \Lambda^2(S,s) } 
 \sqrt {s}  \| \beta_{{\cal N}_0^*}^* \|_2 +
  { \theta (S,s) \over \Lambda^2(S,s) } \| \beta_{S^c}^* \|_1  $$
  $$ = \vartheta_{\rm weak-RIP} (S,s)  \biggl ( \sqrt {s}
    \| \beta_{{\cal N}_0^*}^* \|_2 +  \| \beta_{S^c}^* \|_1 \biggr ) . $$
 Thus, 
   $$
  ( \beta_{S^c}^* )^T
 \Sigma_{2,1} \Sigma_{1,1}^{-1} (S) \tau_S^*  -  \| \beta_{S^c}^* \|_1 $$
 $$ \le\vartheta_{\rm weak-RIP} (S,s)
  \sqrt {s}  \| \beta_{{\cal N}_0^*}^* \|_2 - (1-\vartheta_{\rm weak-RIP} (S,s))
   \| \beta_{S^c}^* \|_1$$ $$ \le \vartheta_{\rm weak-RIP} (S,s)
  \sqrt {s}  \| \beta_{{\cal N}_0^*}^* \|_2.
   $$
    Hence, by Lemma \ref{start.lemma}, 
    $$2 \| ( f_{\beta_{S^c}^*} )^{A_S} \|^2 \le  
     \vartheta_{\rm weak-RIP} (S,s)\biggl (  \lambda   \sqrt {s}  \| \beta_{{\cal N}_0^*}^* \|_2
     \biggr) .$$

 \hfill $\sqcup \mkern -12mu \sqcap$
 
 \begin{lemma} \label{inproduct.lemma} Suppose that
 $$\vartheta_{\rm weak-RIP} (S,s) < 1 . $$ 
 Then for any subset $\tilde {\cal N} \subset S^c$, with $| \tilde {\cal N} | \le s$,
 and any $b \in \R^p$
 $$| (f_{b_{\tilde {\cal N}}} , f^*-f^0) | \le  { \lambda \sqrt {s} \over \phi(S,2s)  \Lambda(S,s)}
 \left ( \theta(S,s) +  \sqrt {(1+ \delta_{s,s})\theta(S,s)/2} \right ) 
 \| b_{\tilde {\cal N}} \|_2 .
 $$
 \end{lemma}
 
 {\bf Proof of Lemma \ref{inproduct.lemma}.} We have
 $$ | (f_{b_{\tilde {\cal N}}} , f^*- f^0) | \le
 |  (f_{b_{\tilde {\cal N}}} , (f^*-f^0)^{P_S}) |+| (f_{b_{\tilde {\cal N}}} , (f^*)^{A_S}) |$$
 Let us write
 $$(f^*-f^0)^{P_S} := f_{\gamma_S}. $$
 Then, invoking Lemma \ref{oracle.lemma}, 
 $$\| \gamma_S \|_2 \le \| f_{\gamma_S} \|/ \Lambda(S,s) =
 \| (f^*-f^0)^{P_S} \|/ \Lambda(S,s) \le \| f^* - f^0 \|/ \Lambda(S,s)
 $$ $$ \le \lambda \sqrt {s} / \biggl (\phi(S,2s)  \Lambda(S,s) \biggr ). $$
 It follows that

\newpage
$$   | (f_{b_{\tilde {\cal N}}} , (f^*-f^0)^{P_S}) | \le \theta(S,s)
\| b_{\tilde {\cal N}} \|_2 \| \gamma_{S} \|_2  $$ $$\le
 \theta(S,s)  \| b_{\tilde {\cal N}} \|_2 \lambda \sqrt {s} / \biggl (\phi(S,2s)  \Lambda(S,s) \biggr )
  .$$
 Moreover, we have
 $$\| \beta_{{\cal N}_0^*}^* \|_2 \le  \| \beta_{{\cal N}_*}^*- \beta_{{\cal N}_*}^0 \|_2 \le
 \lambda \sqrt {s} / \phi^2 (S,2s) . $$
 So, by Lemma \ref{antiprojection.lemma}, 
 $$ \| (f_{\beta_{S^c}^*})^{A_S}\|^2\le  { \theta (S,s) \over \Lambda^2(S,s) } 
\lambda \sqrt {s}  \| \beta_{{\cal N}_*}^* - \beta_{{\cal N}_*}^0 \|_2 /2 $$
 $$ \le \lambda^2 s \theta(S,s) /\biggl  ( 2 \phi^2 (S,2s) \Lambda^2 (S,s)\biggr ) . $$
 Therefore
 $$| (f_{b_{\tilde {\cal N}}} , (f^*)^{A_S}) |\le
 \| f_{b_{\tilde {\cal N}}} \| \| (f^*)^{A_S}) \| 
  \le \lambda \sqrt{ s}
 \sqrt{ \theta(S,s)/2} /\biggl  (  \phi (S,2s) \Lambda (S,s)\biggr )  \| f_{b_{\tilde {\cal N}}} \|$$
 $$ \le { \sqrt {(1+ \delta_{s})\theta(S,s)/2} \over
 \phi (S,2s) \Lambda (S,s) } \lambda \sqrt {s} \| b_{\tilde {\cal N}} \|_2. $$
 
The next result shows that if the constants are small enough,
then there will be no more than $s$ false positives.
We define
\begin{equation}\label{alpha-eq}
\alpha (S):=  { \left ( \sqrt {2} \theta(S,s) +  \sqrt {(1+ \delta_{s})\theta(S,s)} \right ) \over \phi(S,2s) 
\Lambda(S,s) }  .
\end{equation}

\begin{lemma} \label{RIPselect.lemma} 
Suppose that 
$$\alpha(S) < 1 . $$
Then $| S_* \backslash S | < s $.
\end{lemma}

{\bf Proof of Lemma \ref{RIPselect.lemma}}
Since $\alpha(S) < 1$, Lemma
\ref {inproduct.lemma} implies that  for any $\tilde {\cal N}\subset S^c$, with
$| \tilde {\cal N}|\le s $, and for any $b$ with $\| b_{\tilde {\cal N}} \|_2
\not= 0 $,  
$$| ( f_{b_{\tilde {\cal N}}} , f^* - f_0) | 
< \lambda \sqrt {s/2}\| b_{\tilde {\cal N}} \|_2  . $$
Hence, taking $b_j = (\psi_j, f^*- f^0)$, $j \in \tilde {\cal N}$,
$$ \sum_{j \in \tilde {\cal N}} |(\psi_j, f^* - f^0 ) |^2< 
  \lambda^2 s/2  .$$

For $j \in S_*\backslash S $ we have by the KKT conditions
$$|2 ( \psi_j , f^* - f^0 ) | \ge \lambda . $$
Suppose now that $|S_* \backslash S| \ge s$. 
Then there is a subset  ${\cal N}^{\prime}$ of $S_* \backslash S $, with size
$| {\cal N}^{\prime} | = s$, and
we have
$$  \lambda^2 s / 2 >  \sum_{j \in {\cal N}^{\prime} } 
| ( \psi_j , f^* - f^0
 ) |^2 \ge \lambda^2 | {\cal N}^{\prime} |/2. $$
This is a contraction, and hence $| S_*\backslash S | < s $.
\hfill $\sqcup \mkern -12mu \sqcap$

\vskip .1in
This leads to the following result.

\begin{theorem}\label{RIPirr.theorem}
Suppose that $\alpha(S)<1$, see (\ref{alpha-eq}). 
 Then the weak $(S,2s)$-irrepresentable condition holds.
\end{theorem}

{\bf Proof of Theorem \ref{RIPirr.theorem}.} 
As $\alpha(S) <1$, we know that $\phi(S,2s)>0$.
Take an arbitrary
$\tau_S^0 \in \{ -1 , 1 \}^s $, and a
$\beta_0$ satisfying $\beta_S^0 = \beta^0$, ${\rm sign} (\beta_S^0) = \tau_S^0$, and
$$|\beta^0|_{\rm min} > \lambda \sqrt s / \phi^2 (S,2s). $$
By Lemma \ref{oracle.lemma}, the Lasso satisfies
$$ \| \beta_S^* - \beta_S^0 \|_2 \le \lambda \sqrt {s} / \phi^2 (S,2s) . $$
Hence, we must have $S_* \supset S$,
and $\tau_S^* = \tau_S^0$. Moreover, by Lemma \ref{RIPselect.lemma},
$|S_*| < 2s $. By Part 3 of Lemma \ref{select.lemma}, we must have
$$\| \Sigma_{2,1} (S_*) \Sigma_{1,1}^{-1} (S_*) \tau_{S_*}^* \|_{\infty} \le 1 . $$
Since $\tau_{S}^0=\tau_S^*$ is arbitrary and $\tau_{S_*}^* \in \{ -1, 1\}^{|S_*|}$,
we conclude that the weak $(S,2s)$-irrepresentable condition holds
(in fact the weak $(S, 2s-1)$-irrepresentable condition holds). 

\hfill $\sqcup \mkern -12mu \sqcap$

\begin{corollary}\label{RIPirr.corollary}
The RIP is the condition $\vartheta_{\rm RIP}<1$, or equivalently
$$\delta_s + \theta_{s,s} + \theta_{s, 2s} < 1 . $$
\cite{candes2005decoding} show that $\delta_{2s} \le \theta_s + \delta_s$.
The restricted isometry constant $\delta_s$ has to be less than one, so
we may use the bound $1+ \delta_s \le 2 $.
Moreover, it is clear that $\theta (S,N) \le \theta_{s,N}$, and 
$\Lambda^2(S,N) \ge 1- \delta_{N}$.  
Inserting these bounds in Corollary \ref{RIP.corollary}
we find
$$\phi(S,2s) \Lambda(S,s)  \ge (1-\delta_s - \theta_{s,s} - \theta_{s,2s}) 
\sqrt {1-\delta_s \over 1- \delta_s - \theta_{s,s} } \ge  
(1-\delta_s - \theta_{s,s} - \theta_{s,2s}) 
. $$
It follows that
$$\alpha(S) \le  {\sqrt {2}( \theta_{s,s} + \sqrt {\theta_{s,s}) } \over
1- \delta_s - \theta_{s,s} - \theta_{s,2s} } 
 . $$
 For example, if $\delta_s \le \sqrt{2}-1$ and $\theta_{s, 2s} \le {1 \over 16} $, we get
 (invoking $\theta_{s, s} \le \theta_{s, 2s}$)
 $$\alpha(S) \le 0.96. $$
 We conclude that the RIP with small enough constants implies the
 weak $(S,2s)$-irrepresentable condition.
\end{corollary}

As \cite{candes2005decoding} show, the RIP implies exact recovery.
To complete the picture, we now show that the $(S,s)$-irrepresentable condition
also implies exact recovery.

The linear programming problem is
$$\min \{ \| \beta \|_1 : \ \| f_{\beta} - f^0 \|= 0 \} , $$
where, as before $f^0 = f_{\beta^0} $ with $\beta^0 = \beta_S^0 $. 
Let $\beta^{\rm LP}$ be the minimizer of the linear programming problem.

\begin{lemma}\label{recovery.lemma}
Suppose the $(S,s)$-irrepresentable condition holds. Then one has exact recovery, i.e., $\beta^{\rm LP} = \beta^0 $.

\end{lemma}

{\bf Proof of Lemma \ref{recovery.lemma}.} This follows from 
\cite{candes2005decoding}. They show that $\beta^{\rm LP}= \beta^0$ if one
can find a $g \in L_2(P)$, such that\\
(i) $( \psi_j , g )= \tau_j^0 $, for all  $j \in S $,\\
(ii) $| (\psi_j , g) |< 1 $ for all $j \notin S $,\\
where, as before, $\tau_S^0 := {\rm sign} (\beta_S^0)$.
The $(S,s)$-irrepresentable condition says that this is true for
$g= f_{b_S} $, where $b_S = \Sigma_{1,1}^{-1} (S) \tau_S^0 $.
\hfill $\sqcup \mkern -12mu \sqcap$

\section{ The $(S,s)$-uniform irrepresentable condition implies the $S$-compatibility condition}
\label{irrcomp.section}
As the $(S,s)$-irrepresentable condition implies variable
selection, one expects it will be more restrictive than the
compatibility condition, which only implies a bound
for the prediction error (and $\ell_1$-estimation error). 
This turns out to be indeed the case,
albeit we prove it only under the uniform version of the irrepresentable
condition. 

\begin{theorem} \label{irrcomp.theorem}
Suppose that
$$\vartheta_{\rm irrepresentable} (S,s) < 1/L . $$
Then
$$ \phi_{\rm compatible}^2 (L,S) 
\ge (1-L \vartheta_{\rm irrepresentable}(S,s))^2 \Lambda^2 ( S,s) .$$
\end{theorem}

{\bf Proof of Theorem \ref{irrcomp.theorem}.} 
Define,
$$\beta^{\diamond} := \arg \min_{\beta} \{ \| f_{\beta} \|^2 :  \ \| \beta_S \|_1 =1 , \| \beta_{S^c} \|_1 \le L \}. $$
Let us  write
$f^{\diamond}: = f_{\beta}^{\diamond}$, $f_S^{\diamond} := f_{\beta_S}^{\diamond}$ and
$f_{S^c}^{\diamond} := f_{\beta_{{S^c}}^{\diamond}} $.
Introduce a Lagrange multiplier $\lambda \in \R$ for the constraint
$\|\beta_s\|_1 = 1$. By the KKT conditions,
there exists a vector $\tau_S^{\diamond}$, with $\| \tau_{S}^{\diamond} \|_{\infty} \le 1 $,
such that $\tau_S^T \beta_S^{\diamond}= \| \beta_S^{\diamond} \|_1$, and such that
\begin{equation}\label{KKT-proof9.1}
\Sigma_{1,1} (S) \beta_S^{\diamond} +  \Sigma_{1,2}(S)
\beta_{S^c}^{\diamond} = - \lambda \tau_S^{\diamond} .
\end{equation}

 By multiplying by $(\beta_S^{\diamond})^T$, we obtain
 $$ \| f_S^{\diamond} \|^2 + ( f_S^{\diamond} , f_{S^c}^{\diamond}) 
 = - \lambda \| \beta_S^{\diamond} \|_1 . $$
The restriction $\| \beta_S^{\diamond} \|_1=1$ gives
 $$ \| f_S^{\diamond}\|^2 + ( f_S^{\diamond} , f_{S^c}^{\diamond} ) 
 = - \lambda  . $$
We also have from (\ref{KKT-proof9.1})
\begin{equation}\label{inverse.eq}
\beta_S^{\diamond}+ \Sigma_{1,1}^{-1} (S) \Sigma_{1,2} (S) \beta_{S^c}^{\diamond} =-\lambda 
\Sigma_{1,1}^{-1} \tau_S^{\diamond} . 
\end{equation}
Hence, by multiplying with $(\tau_S^{\diamond})^T$,
$$\| \beta_S^{\diamond} \|_1 +( \tau_S^{\diamond})^T  \Sigma_{1,1}^{-1} (S) \Sigma_{1,2}(S) \beta_{S^c}^{\diamond} =
- \lambda (\tau_S^{\diamond})^T \Sigma_{1,1}^{-1} \tau_S^{\diamond} , $$
or
$$1 =-  (\tau_S^{\diamond})^T  \Sigma_{1,1}^{-1} (S) \Sigma_{1,2} (S) \beta_{S^c}^{\diamond}-
\lambda (\tau_S^{\diamond})^T \Sigma_{1,1}^{-1}(S) \tau_S^{\diamond} $$
$$ \le \vartheta \| \beta_{S^c}^{\diamond} \|_1 -
\lambda (\tau_S^{\diamond})^T \Sigma_{1,1}^{-1}(S) \tau_S^{\diamond} $$
$$\le L  \vartheta  -
\lambda (\tau_S^{\diamond})^T \Sigma_{1,1}^{-1} (S) \tau_S^{\diamond} .$$
Here, we applied that  the $(S,s)$-uniform irrepresentable condition, 
with $\vartheta = \vartheta_{\rm irrepresentable}(S,s)$, and
the condition $\| \beta_{S^c} \|_1 \le L $.
Thus
$$1-L \vartheta   \le - \lambda  (\tau_S^{\diamond})^T \Sigma_{1,1}^{-1} (S) \tau_S^{\diamond} . $$
Because $1-L \vartheta  >0$ and $(\tau_S^{\diamond})^T \Sigma_{1,1}^{-1} (S) 
\tau_S^{\diamond} \ge 0$, 
this implies that $\lambda < 0 $, and in fact that
$$(1-L \vartheta )  \le -\lambda s / \Lambda^2 (S,s),$$
where we invoked
$$( \tau_S^{\diamond})^T \Sigma_{1,1}^{-1} (S) \tau_S^{\diamond} \le \| \tau_S^{\diamond} \|_2^2
/ \Lambda^2 (S,s) 
\le s/  \Lambda^2 (S,s) . $$
So
$$- \lambda \ge (1-L \vartheta) \Lambda^2 (S,s) / s .$$

Continuing with (\ref{inverse.eq}), we moreover have
$$(\beta_{S^c}^{\diamond} )^T \Sigma_{2,1}(S) \beta_S^{\diamond} +
(\beta_{S^c}^{\diamond} )^T \Sigma_{2,1} (S) \Sigma_{1,1}^{-1} (S) \Sigma_{1,2}(S)
\beta_{S^c}^{\diamond}$$ $$ =-\lambda  (\beta_{S^c}^{\diamond})^T \Sigma_{2,1} (S)
\Sigma_{1,1}^{-1} (S) \tau_S^{\diamond} . $$
In other words,
$$( f_S^{\diamond}, f_{S^c}^{\diamond} ) + \| (f_{S^c}^{\diamond} )^{P_S} \|^2 
=-\lambda  (\beta_{S^c}^{\diamond})^T \Sigma_{2,1} (S)
\Sigma_{1,1}^{-1} (S) \tau_S^{\diamond} , $$
where  $( f_{S^c}^{\diamond} )^{P_S} $ is the projection of
$f_{S^c}^{\diamond}$ on the space spanned by $\{ \psi_k \}_{k \in S} $.
Again, by the $(S,s)$-uniform irrepresentable condition and by
$\| \beta_{S^c}^{\diamond} \|_1 \le L$, 
$$ \left | (\beta_{S^c}^{\diamond})^T \Sigma_{2,1} (S)
\Sigma_{1,1}^{-1} (S) \tau_S^{\diamond}\right | \le  \vartheta \|\beta_{S^c}^{\diamond} \|_1 \le
L\vartheta  , $$
so 
$$-\lambda  (\beta_{S^c}^{\diamond})^T \Sigma_{2,1} (S)
\Sigma_{1,1}^{-1} (S) \tau_S^{\diamond} 
=|\lambda |  (\beta_{S^c}^{\diamond})^T \Sigma_{2,1} (S)
\Sigma_{1,1}^{-1} (S) \tau_S^{\diamond} 
$$ $$ \ge- |\lambda |\left |   (\beta_{S^c}^{\diamond})^T \Sigma_{2,1} (S)
\Sigma_{1,1}^{-1} (S) \tau_S^{\diamond} \right |
 \ge- | \lambda|  L \vartheta  = \lambda L \vartheta  .  $$
It follows that
$$ \| f^{\diamond} \|^2 = \| f_{S}^{\diamond} \|^2 + 2 (f_S^{\diamond}, f_{S^c}^{\diamond}) + \| f_{S^c}^{\diamond} \|^2 $$
$$ = - \lambda +(f_S^{\diamond}, f_{S^c}^{\diamond}) + \| f_{S^c}^{\diamond}\|^2$$
$$ \ge -\lambda + (f_S^{\diamond}, f_{S^c}^{\diamond}) + \|( f_{S^c}^{\diamond})_S^{\rm P} \|^2
\ge -\lambda + \lambda L\vartheta  =- \lambda (1-L \vartheta )  $$ $$
\ge (1-L\vartheta  )^2  \Lambda^2 ( S,s ) / s .$$
Finally note that $\| f^{\diamond} \|^2 = \phi_{\rm compatible}^2(L,S)/s$. 
\hfill $\sqcup \mkern -12mu \sqcap$

\section{Verifying the compatibility and restricted eigenvalue
  condition} \label{verify.section} 
%
%
%
%
 
In this section, we discuss the theoretical verification of the conditions.
Determining a restricted $\ell_1$-eigenvalue is in itself again
a Lasso type of problem. Therefore, it is very useful
to look for some good lower bounds.

A first, rather trivial, observation is that if 
$\Sigma$ is non-singular, the restricted eigenvalue condition holds for all
$L$, $S$ and $N$, with $\phi^2(L,S,N) \ge \Lambda_{\rm min}^2 (\Sigma)$, the
latter being the smallest eigenvalue of $\Sigma$. If $\Sigma$ is the
population covariance matrix of a random design, i.e., the probability measure $Q$
is the theoretical distribution of observed co-variables in ${\cal X}$, assuming
positive definiteness of $\Sigma$ is not very restrictive. We will present
some examples in Section \ref{matrices.section}. 
Compatibility conditions for the population Gram matrix are of direct
relevance if one replaces $L_2$-loss by robust convex loss \citep{geer08}.
But, as we will show in the next subsection, even if $\Sigma$
corresponds to the empirical 
covariance matrix of a fixed design, i.e., the measure $Q$ is the empirical
measure $Q_n$ of $n$ observed co-variables in ${\cal X}$, the compatibility and
restricted eigenvalue condition is often ``inherited" from the
population version.
Therefore, even for
fixed designs (and singular $\Sigma$), the collection  of cases where 
compatibility or restricted eigenvalue conditions hold is quite large.

\subsection{Approximating the Gram matrix}\label{subsec.appgram}

For two (positive semi-definite) matrices $\Sigma_0$ and
$\Sigma_1$, we define the supremum distance
$$d_{\infty} ( \Sigma_1,  \Sigma_0)  := \max_{j,k} | ( \Sigma_1)_{j,k} -( \Sigma_0)_{j,k} |.$$

Generally, perturbing  the entries in $\Sigma$ by a small amount
may have a large impact on the eigenvalues of $\Sigma$. This is not true for
(adaptive) restricted $\ell_1$-eigenvalues, as is shown in the next lemma and
its corollary. 

\vskip .1in
\begin{lemma} \label{metrics} 
Assume
$$d_{\infty} ( \Sigma_1 , \Sigma_0 ) \le \tilde \lambda . $$
Then $\forall \ \beta \in {\cal R} (L,S)$, 
$$ \left | { \| f_{\beta} \|_{ \Sigma_1}^2 \over \|  f_{\beta} \|_{\Sigma_0}^2}  -1 \right | \le 
{(L+1)^2 \tilde \lambda s \over  \phi_{\rm compatible}^2 (\Sigma_0 , L,S)}, 
 $$
and similarly, 
 $\forall \  {\cal N} \supset S, \
|{\cal N}| = N $, and $\forall \ \beta \in {\cal R}(L,S,{\cal N}) $, 
$$
 \left | { \| f_{\beta} \|_{ \Sigma_1}^2 \over \|  f_{\beta} \|_{\Sigma_0}^2}  -1 \right | \le 
{(L+1)^2 \tilde \lambda s \over  \phi^2 (\Sigma_0 , L,S,N)}, 
 $$
and
$\forall \  {\cal N} \supset S, \
|{\cal N}| = N $, and $\forall \ \beta \in {\cal R}_{\rm adaptive} (L,S,{\cal N}) $, 
$$
 \left | { \| f_{\beta} \|_{ \Sigma_1}^2 \over \|  f_{\beta} \|_{\Sigma_0}^2}  -1 \right | \le 
{(L+1)^2 \tilde \lambda s \over  \phi_{\rm adaptive}^2 (\Sigma_0 , L,S,N)}
 .$$
\end{lemma}

{\bf Proof of Lemma \ref{metrics}.} 
For all $\beta$,
$$ \left | \| f_{\beta} \|_{\Sigma_1}^2 -\| f_{\beta} \|_{\Sigma_0}^2 \right | =
| \beta^T  \Sigma_1 \beta - \beta^T \Sigma_0 \beta | $$
$$= | \beta^T (  \Sigma_1 - \Sigma_0 ) \beta | \le
\tilde \lambda \| \beta \|_1^2 . $$
But if $\beta \in {\cal R} (L,S)$, it holds that
$\| \beta_{S^c} \|_1 \le L \| \beta_S \|_1 $, and hence
$$\| \beta \|_1 \le (L+1) \| \beta_{ S} \|_1 \le 
(L+1) \| f_{\beta} \|_{\Sigma_0} \sqrt {s}  /\phi_{\rm compatible} (\Sigma_0, L,S) . $$
This gives
$$ \left | \| f_{\beta} \|_{\Sigma_1}^2 -\| f_{\beta} \|_{\Sigma_0}^2 \right | \le
(L+1)^2  \tilde \lambda \| f_{\beta } \|_{\Sigma_0}^2 s /\phi_{\rm compatible}^2 (\Sigma_0 , L , S) . $$
The second result can be shown in the same way,
and the third result as well as for
$\beta \in {\cal R}_{\rm adaptive} ( L, S, {\cal N})$, it holds that
$\| \beta_{S^c} \|_1 \le L \sqrt s \| \beta_S \|_2 $, and hence
$$\| \beta \|_1 \le L \sqrt s \| \beta_S \|_2 + \| \beta_{S }\|_1 \le
(L+1) \sqrt s \| \beta_S \|_2 . $$
\hfill $\sqcup \mkern -12mu \sqcap $

\begin{corollary}\label{metrics2} 
We have
$$\phi_{\rm compatible}(\Sigma_1, L, S) \ge
\phi_{\rm compatible} (\Sigma_0 , L, S) - (L+1) \sqrt {d_{\infty}(\Sigma_0 , \Sigma_1) s } . $$
Similarly, 
$$\phi(\Sigma_1, L, S,N) \ge
\phi(\Sigma_0 , L, S,N) - (L+1) \sqrt {d_{\infty}(\Sigma_0 , \Sigma_1) s } , $$
and the same result holds for the adaptive version.
\end{corollary}

Corollary \ref{metrics2} shows that
if one can find a matrix $\Sigma_0$  with well-behaved smallest eigenvalue, in a small enough $\ell_{\infty}$-neighborhood of $\Sigma_1$,
then the restricted eigenvalue condition
holds for $\Sigma_1$. As an example, 
consider  the situation where $\psi_j(x)
= x_j\ (j=1,\ldots ,p)$ and where 
$$ \hat \Sigma := {\bf X}^T {\bf X} / n= (\hat \sigma_{j,k}) , $$
where ${\bf X}= ( X_{i,j} )$ is a $(n \times p)$-matrix whose columns consist
of i.i.d.\ ${\cal N}(0,1)$-distributed entries (but allowing for dependence
between columns). We denote by $\Sigma$ the population covariance matrix of a
row of ${\bf X}$. 
Using a union bound, it is not difficult to show that for
all $t >0$,  and for 
$$\tilde \lambda (t)  := \sqrt { 4 t + 8 \log p  \over n} + { 4 t + 8 \log p  \over n}, $$
one has the inequality
\begin{equation}\label{gauss-randomd}
\PP \biggl  ( d_{\infty}( \hat \Sigma , \Sigma) \ge
\tilde \lambda (t)  \biggr ) \le 2 \exp[-t] .
\end{equation}
This implies that if the smallest eigenvalue $\Lambda_{\rm min}^2(\Sigma)$ of
$\Sigma$ is bounded away from zero, and if the
sparsity $s$ is of smaller order $o(\sqrt {n / \log p }) $,
then the restricted eigenvalue condition holds with
constant $\phi(S,N)$ not much smaller than
$\Lambda_{\rm min}(\Sigma)$. 
The result can be extended to distributions with
Gaussian tails.

\subsection{Some examples}\label{matrices.section}
In the following, our discussion mainly applies for $\Sigma$ being the
population covariance matrix. For $\Sigma$ being the empirical covariance matrix, the
assumptions in the discussion below are unrealistic, but as seen in the previous
section, the population properties can have
important implications for  the restricted eigenvalues of the empirical covariance matrix.

\begin{example}\label{ex.equicorr}
Consider the matrix
$$\Sigma := ( 1-\rho )I + \rho \iota \iota^T , $$
with $0< \rho < 1$, and $\iota := ( 1 , \ldots , 1)^T$ a vector of 1's.
Then the smallest eigenvalue of $\Sigma$ is $\Lambda_{\rm min}^2 (\Sigma) =1-
\rho$, so the $(L,S,N)$-restricted eigenvalue condition
holds with $\phi^2(L,S,N) \ge 1-\rho$. The uniform $(S,s)$-irrepresentable condition is always met.
The largest eigenvalue of
$\Sigma$ is $(1-\rho) + \rho p$.  Hence, the restricted isometry constants
$\delta_s$ are defined only for $\rho < 1/(s-1)$.
\end{example}

\begin{example}\label{ex.toeplitz}
In this example, $\Sigma$ is a Toeplitz matrix, defined as follows. Consider
a positive definite function 
$$ R(k),\ k \in \Z ,$$
which is symmetric ($R(k) = R(-k)$) and sufficiently regular in the following
sense. The corresponding spectral density
$$f_{\rm spec}(\gamma)  := \sum_{k=-\infty}^{\infty} R(k) \exp(-ik \gamma)\
(\gamma \in [-\pi,\pi])$$
is assumed to exist, to be continuous and periodic, and
$$ \gamma_0 := \argmin_{\gamma \in [0,\pi]} f_{\rm spec}(\gamma)$$
is assumed unique, with $f(\gamma_0) = M > 0$. Moreover, we suppose that $f_{\rm
  spec}(\cdot)$ is $(2 \alpha)$ continuously differentiable at $\gamma_0$,
with $f^{(2 \alpha)}(\gamma_0) > 0$. 
A Toeplitz matrix is
$$\Sigma =  (\sigma_{j,k}),\ \sigma_{j,k} := R(|j-k|),$$
where $R(\cdot)$ satisfies the conditions described above (in terms of the
spectral density). A special case arises with $\sigma_{j,k} = \rho^{|j-k|}$
for some $0 \le \rho < 1$.
The smallest eigenvalue $\Lambda_{\rm min}^2 (\Sigma)$ of $\Sigma$ is bounded away
from zero where the bound is independent of $p$ \citep{Parter:61}.


\end{example}

\begin{example}\label{ex.block}
Consider a matrix $\Sigma$ which is of block structure form:
$$\Sigma =\ \mbox{{\rm diag}}(\Sigma_1,\ldots, \Sigma_k),$$
where the $\Sigma_j$ are $(m \times m)$ covariance matrices ($j=1,\ldots ,k$) 
(the restriction to having the same dimension $m$ can be easily dropped) and
$km = p$. If the minimal eigenvalues satisfy
$$\min_j\Lambda_{\rm min}^2(\Sigma_j) \ge \eta^2 > 0,$$
then the minimal eigenvalue of $\Sigma$ is also bounded from below by $\eta^2 >
0$. When $m$ is much smaller than $p$, it is (much) less restrictive that small
$m \times m$ covariance matrices $\Sigma_j$ have well-behaved minimal
eigenvalues than large $p \times p$ matrices.

\end{example}

\begin{example} \label{irr.example} 
This example presents a case where the compatibility condition holds, but where
the uniform irrepresentable constant is
very large. We also calculate the adaptive restricted regression.
Let the first $s$ indices $ \{ 1 , \ldots , s\}$ be the active set $S$ and
suppose that 
 $$ \Sigma := \pmatrix { I &  \Sigma_{1,2} \cr  \Sigma_{2,1} & \Sigma_{2,2} \cr},
$$
where  $I$ is the $(s \times s)$-identity matrix,
and
$$ \Sigma_{2,1}:= \rho( b_2 b_1^T)  , $$
with $0\le  \rho < 1$, and with $b_1$ an $s$-vector and $b_2$ a $(p-s)$-vector,
satisfying $\| b_1 \|_2=\| b_2 \|_2=1 $.
Moreover, 
$ \Sigma_{2,2} $ is some $(p-s)\times (p-s) $-matrix,
with ${\rm diag} ( \Sigma_{2,2} ) = I$, and 
with smallest eigenvalue $ \Lambda_{\rm min}^2 (\Sigma_{2,2}) $.
One easily verifies that
$$\Lambda_{\rm min}^2 (\Sigma) \ge  \Lambda_{\rm min}^2 (  \Sigma_{2,2} ) - \rho. $$
Moreover, for
$b_1 := ( 1 , 1, \ldots , 1)^T/\sqrt {s}$
and 
$b_2:= (1 , 0 , \ldots , 0 )^T$, and $\rho > 1 / \sqrt {s} $, 
the $(S,s)$-uniform irrepresentable condition does not hold,
as in that case
$$
\sup_{\| \tau_S \|_{\infty} \le 1 }
\| \Sigma_{2,1} (S) \Sigma_{1,1}^{-1} (S) \tau_S \|_{\infty} = \rho \sqrt {s} .$$
However, for any $N >s$, the $(S,N)$-uniform irrepresentable condition does hold.
We moreover have
$$\vartheta_{{\rm adaptive}} (S)=
{  \sqrt {s} \| \Sigma_{1,2} \|_{2, \infty} =\sqrt s \rho 
  } ,$$
 i.e. (since $\Lambda (S,s)=1$), the bounds of Lemma \ref{coherencelemma}
 and Theorem \ref{thetabound.theorem} are strict in this example.
  
\end{example}

\begin{example}\label{ex.compatib}

We recall that $\phi_{\rm compatible} (S) \ge
\phi(S,s)$. 
Here is an example where the compatibility
condition holds with reasonable $\phi_{\rm compatible}^2 (S)$, but where the
restricted eigenvalue $\phi^2(S,s)$ is very small. Assume $s > 2$.
Let the first $s$ indices $ \{ 1 , \ldots , s\}$ be the active set $S$ with
corresponding $(s \times s)$ covariance matrix $\Sigma_{1,1}$, and
suppose that 
 $$ \Sigma := {\rm diag} (\Sigma_{1,1} , I ) , 
$$
where
$$ \Sigma_{1,1} = {\rm diag} (B,I) , $$
and, for some $0 \le  \rho < 1- 1/(s-2)$, 
$$B = \pmatrix { 1 & \rho \cr \rho & 1 \cr } . $$
We then have
$$\beta_S^T \Sigma_{1,1} \beta_S =
(1- \rho ) ( \beta_1^2 + \beta_2^2) + \rho ( \beta_1 + \beta_2)^2 +
\sum_{j=3}^s \beta_j^2 $$ $$  \ge
(1- \rho ) ( \beta_1^2 + \beta_2^2) + 
(\sum_{j=3}^s |\beta_j|)^2  /(s-2)
$$
Hence, 
$$\min_{\| \beta_S \|_1 =1 } \beta_S^T \Sigma_{1,1} \beta_S \ge
\min_{|\beta_1| + |\beta_2 | \le 1 } \biggl \{ 
(1- \rho ) ( \beta_1^2 + \beta_2^2) + (1- | \beta_1| -| \beta_2| )^2 /(s-2)  \biggr \}
  $$
  $$\ge \min_{|\beta_1| + |\beta_2 | \le 1 } \biggl  \{ \sum_{j=1,2} \biggl (1- \rho+ {1 \over s-2} 
  \biggr ) \beta_j^2 +
 {1 \over 2-s} -
  2 {|\beta_1| + |\beta_2 | \over s-2} \biggr \} $$
  $$= \min_{|\beta_1| + |\beta_2 | \le 1 } \biggl \{ 
   { (s-2 )(1-\rho) +1 \over s-2} \sum_{j=1,2} \biggl ( | \beta_j|- {1 \over (s-2)(1-\rho)+1 }\biggr  )^2 
     \biggr \} $$
 $$- {2 \over (s-2) \biggl  ((s-2)(1-\rho) +1\biggr )} +{1 \over s-2} $$
  $$\ge { (s-2) (1-\rho) -1 \over (s-2) \biggl ((s-2)(1-\rho) +1\biggr ) } . $$
  It follows that
  $$\phi_{\rm compatible}^2 (S) =
  \min_{\| \beta_S\|_1 = 1 ,\  \| \beta_{S^c} \|_1 \le 1 }
  { s \beta^T \Sigma \beta \over \| \beta_S \|_1^2 } \ge
   { s\biggl ((s-2) (1-\rho) -1 \biggr ) \over (s-2) \biggl ((s-2)(1-\rho) +1\biggr ) } $$
   $$\ge { (s-2) (1-\rho) -1  \over(s-2)(1-\rho) +1 } . $$
  On the other hand
  $$\phi^2(S,s) = \Lambda^2 (S,s) = (1-\rho) . $$
  Hence, for example when $1-\rho=3/(s-2)$, we get
  $$\phi_{\rm compatible}^2  (S) \ge 1/2$$
  and
  $$\phi^2(S,s) = {3 \over s-2} . $$
  Clearly, for large $s$, this means that $\phi_{\rm compatible} (S)$
  is much better behaved than $\phi(S,s)$. Note that large $s$ in this example
  (with $1-\rho=3/(s-2)$) corresponds to a correlation $\rho$ close to one,
  i.e., to a case where $\Sigma$ is ``almost" singular. 
  
 \end{example}

\section{Adding noise}\label{noise.section}
We now consider the Lasso estimator based on $n$ noisy observations.
Let $X_i \in {\cal X}$ ($i=1 , \ldots , n $) be the co-variables,
and $Y_i \in \R$ ($i=1 , \ldots , n $) be the response variables.
The noisy Lasso is 
$$\hat\beta := \arg \min_{\beta} \biggl \{ {1 \over n}
\sum_{i=1}^n | Y_i - f_{\beta} (X_i) |^2   + \lambda \| \beta \|_1 \biggr \} . $$
The design matrix is
$${\bf X} = {\bf X}_{n \times p } : = 
( \psi_j (X_i)). $$
The empirical Gram matrix is
$$\hat \Sigma := {\bf X}^T {\bf X} / n = \int \psi^T \psi d Q_n =
( \hat \sigma_{j,k} ), $$
where $Q_n$ is the empirical measure $Q_n := \sum_{i=1}^n \delta_{X_i} /n$.
The $L_2 (Q_n)$-norm is denoted by $\| \cdot \|_n$.
We moreover let $( \cdot , \cdot )_n$ be the $L_2 (Q_n)$-inner product.

 As before, we write $f^0 = f_{\beta^0}$ and now,
$\hat f = f_{\hat \beta} $. 
We consider
$$\epsilon_i := Y_i  - { f}^0 (X_i) , \ i=1 , \ldots , n , $$
as the {\it noise}. Moreover, we write (with some abuse of notation)
$$(f, \epsilon)_n := {1 \over n} \sum_{i=1}^n f (X_i) \epsilon_i , $$
and we define
$$\lambda_0 := 2 \max_{1 \le j \le p }| ( \psi_j , \epsilon )_n |. $$

Here is a simple example which shows how $\lambda_0$ behaves in the case of
i.i.d.\ standard normal errors. 

\begin{lemma} \label{errors.lemma} Suppose that $\epsilon_1 , \ldots , 
\epsilon_n$ are i.i.d. ${\cal N} (0,1)$-distributed, and that
$\hat \sigma_{j,j}  =1$ for all $j$.
Then we have for all $t>0$, and for 
$$\lambda_0 (t) := 2 \sqrt {2t + 2 \log p  \over n } , $$
$$\PP \biggl (2 \max_{1 \le j \le p } | ( \psi_j , \epsilon) _n |  \le  \lambda_0 (t) \biggr   ) \ge 1- 2 \exp [-t] . $$
\end{lemma}

{\bf Proof.} As $\hat \sigma_{j,j}=1$, we know that
$V_j := \sqrt n (  \psi_j , \epsilon)_n $ is ${\cal N}(0,1) $-distributed.
So
$$\PP \left ( \max_{1 \le j \le p} | V_j | > \sqrt { 2t + 2 \log p } \right ) \le
2p \exp \left [ - {2t + 2 \log p \over 2 } \right ] = 2 \exp \left
[ - { t  } \right ] . $$
\hfill $\sqcup \mkern -12mu \sqcap$

\subsection{Prediction error in the noisy case}

 A noisy counterpart of Lemma \ref{oracle.lemma} is:
 
 \begin{lemma} \label{noiseoracle.lemma} Take $\lambda > \lambda_0$, and
 define $L:= ( \lambda + \lambda_0 )/ ( \lambda - \lambda_0) $.
 Then
 $$\| \hat f - f^0 \|_n^2 + {2 \lambda_0 \over L-1} \| \hat \beta_{S^c} \|_1 \le
 { 4 ( L+1)^2 \lambda_0^2 s \over (L-1)^2 \phi_{\rm compatible}^2 (\hat \Sigma , L,S)} . $$
 \end{lemma}
 
 {\bf Proof of Lemma \ref{noiseoracle.lemma}.}
 Because
 $$2 |( \epsilon , \hat f - f^0 ) | \le \biggl ( 2 \max_{1 \le j \le p } | ( \psi_j , \epsilon) | \biggr )
 \| \hat \beta - \beta^0 \|_1 \le \lambda_0  \| \hat \beta - \beta^0 \|_1, $$
 we now have the Basic Inequality
 $$\| \hat f - f^0 \|_n^2 + \lambda \| \hat \beta_{S^c} \|_1 \le 
 \lambda_0 \| \hat \beta - \beta^0 \|_1 + \lambda \| \beta^0 \|_1 . $$
 Hence,
 $$\| \hat f - f^0 \|_n^2 + ( \lambda - \lambda_0) \| \hat \beta_{S^c} \|_1 \le
 (\lambda+ \lambda^0 ) \| \hat \beta_S - \beta_S^0 \|_1 . $$
 Thus,
 $$\| \hat \beta_{S^c} \|_1 \le L \| \hat \beta_S - \beta_S^0 \|_1 . $$
 This implies
 $$ \| \hat \beta_S - \beta_S^0 \|_1 \le \sqrt s \| \hat f - f^0 \|_n /
  \phi_{\rm compatible} (\hat \Sigma , L,S). $$
  So we arrive at
  $$\| \hat f - f^0 \|_n^2 + ( \lambda - \lambda_0) \|\hat  \beta_{S^c} \|_1 \le
 (\lambda+ \lambda^0 ) \sqrt s \| \hat f - f^0 \|_n /
  \phi_{\rm compatible} (\hat \Sigma , L,S). $$
  Now, insert
  $\lambda=\lambda_0 { (L+1)  /( L-1)} $.
\hfill $\sqcup \mkern -12mu \sqcap$
 
 In a similar way, but using $(S,2s)$-restricted eigenvalue conditions, one may prove $\ell_2$-convergence in the noisy case.
 
 Observe that the $S$-compatibility condition now involves the matrix $\hat \Sigma$,
 which is definitely singular when $p > n$. However, we have seen in the previous section
 that, also for such $\hat \Sigma$, compatibility conditions and restricted eigenvalue
 conditions hold in fairly general situations.
 
 \subsection{Noisy KKT}
 
 The KKT conditions in the noisy case become
 $$2 ( \psi_j , \hat f - f^0 )_n -2 ( \psi_j , \epsilon)_n = - \lambda \hat \tau_j, \ 
 j=1 , \ldots , p , $$
 or in matrix notation,
 $$2 \hat \Sigma ( \hat \beta - \beta^0 ) - {\bf X}^T \epsilon /n = - \lambda \hat \tau, $$
 where $\| \hat \tau \|_{\infty} \le 1 $, and $\hat \tau_j  := {\rm sign} (\hat \beta_j )$
 whenever $\hat \beta_j \not= 0 $.
 
To avoid too many repetitions, let us only formulate the noisy version of 
a part of Part 1 of Lemma
\ref{select.lemma}.
 
 \begin{lemma} \label{noisyselect.lemma}
Take $\lambda > \lambda_0$, and
 define $L:= ( \lambda + \lambda_0 )/ ( \lambda - \lambda_0) $.
 Suppose the uniform $(\hat \Sigma, L,S,s)$-irrepresentable condition holds.
 Then $\hat S \subset S$.
 \end{lemma}
 
 {\bf Proof of Lemma \ref{noisyselect.lemma}.}
 This follows from a straightforward  generalization of Lemma \ref{start.lemma},
 where the equalities now become inequalities:
 $$2 \| ( f_{\hat \beta_{S^c}})^{\hat A_S} \|_n^2 \le 
 {2 L \over L-1} \lambda_0 \hat \Sigma_{2,1}(S)  \hat \Sigma_{1,1}^{-1}  (S)\hat \tau_{S} -
 {2 \over L-1}\lambda_0  \| \hat \beta_{S^c} \|_1 . 
 $$
 Here, $f^{\hat A_S}$ is the anti-projection of $f$, in $L_2(Q_n)$, on
 the space spanned by $\{ \psi_j \}_{j \in S}$. 
  
 \hfill $ \sqcup \mkern -12mu \sqcap$
 
 The noisy KKT conditions involve the matrix $\hat \Sigma$.
 Again, as discussed in Subsection \ref{subsec.appgram}, we may replace it by an approximation.
As a consequence,  if this approximation is good enough, we can replace
 $(\hat \Sigma, L, S, s)$-irrepresentable conditions by 
 $(\Sigma, \tilde L , S,s)$-irrepresentable conditions, provided we take
 $\tilde L > L$ large enough.
 
 \begin{lemma} Take $\lambda > \lambda_0$, and
 define $L:= ( \lambda + \lambda_0 )/ ( \lambda - \lambda_0) $.
 Suppose that  \label{approx.lemma}
 $$d_{\infty} ( \hat \Sigma , \Sigma) \le \tilde \lambda , $$
 and
 $$\phi_{\rm compatible} ( \Sigma , L, S) > (L+1) \sqrt {\tilde \lambda s}   $$
 and in fact, that
 $$
 {  (L+1)\sqrt { \tilde \lambda  s}  \over
    \phi_{\rm compatible} (\Sigma , L,S)- (L+1) \sqrt {\tilde \lambda s}} < 1 . $$
   Then
   $$\| (\hat \Sigma - \Sigma) (\hat \beta - \beta^0) \|_{\infty} <
   {2 \lambda_0 \over L-1 } . $$
 \end {lemma} 
 
 {\bf Proof of Lemma \ref{approx.lemma}.}
 We have
 $$\| ( \hat \Sigma - \Sigma )(\hat \beta - \beta^0 ) \|_{\infty} \le
 \tilde \lambda \| \hat \beta - \beta^0 \|_1 \le
 (L+1) \tilde \lambda \| \hat \beta_S - \beta_S^0 \|_1  $$ $$\le
 (L+1)\tilde \lambda  \sqrt {s}  \| \hat f - f^0 \|_n / \phi_{\rm compatible} ( \hat \Sigma , L , S)$$
  $$ \le {2\lambda_0 (L+1)^2 \tilde \lambda  s  \over
   (L-1) \phi_{\rm compatible}^2 (\hat \Sigma , L,S)}  $$
   $$ \le { 2 \lambda_0 (L+1)^2 \tilde \lambda  s \lambda_0 \over
   (L-1) \biggr( \phi_{\rm compatible} (\Sigma , L,S)- (L+1) \sqrt {\tilde \lambda s}\biggl )^2}  .$$
   \hfill $\sqcup \mkern -12mu \sqcap$ 
   
   We conclude that the KKT conditions in the noisy case can be
   exploited in the same way as in the case without noise, albeit that
   one needs to adjust the constants (making the conditions more restrictive).
   
   \section{Discussion} \label{discussion.section}

 We show how various conditions for Lasso oracle results relate to each
 other, as illustrated in Figure \ref{fig1}. Thereby, we also introduce the
 restricted regression condition. 

For deriving oracle results
 for prediction and estimation, the compatibility condition is the
 weakest. Looking at the derivation of the oracle result in Lemma
 \ref{oracle.lemma}, no
 substantial room seems to be left to improve the condition. The restricted
 eigenvalue condition is slightly stronger but in some cases, as
 demonstrated in Example \ref{ex.compatib}, the compatibility condition is
 a real improvement. 
 
For variable selection with the Lasso, the irrepresentable condition is
sufficient (assuming sufficiently large non-zero regression coefficients)
and essentially necessary. We present the, perhaps not unexpected, but
as yet not formally shown, result
that the irrepresentable condition is always stronger than the 
compatibility condition.

We illustrate in Section \ref{verify.section} how - in theory - one can verify the 
compatibility condition. If the sparsity is of small order
$o(\sqrt{n/ \log p })$, we can approximate the empirical Gram matrix by the
population analogue. It is then much more easy and realistic that the
population Gram matrix has sufficiently regular behavior, as illustrated
with our examples in Section \ref{matrices.section}. 
We believe moreover that a sparsity bound of small order
$o(\sqrt{n/ \log p })$ covers a large area of interesting statistical problems.
With larger $s$, the statistical situation is comparable
to one of a nonparametric model with ``(effective) smoothness less than 1/2", 
leading to very slow convergence rates.
In contrast, for example in decoding problems, sparseness up to
the linear-in-$n$ regime can be very important. 
Moreover, in the case of robust convex loss, one
may apply the compatibility condition directly to the
population matrix, i.e., the sparsity regime
$s = o(\sqrt{n/ \log p})$ can be relaxed for such loss functions
(see \cite{geer08}).
We therefore conclude that oracle results for the Lasso
hold under quite general design conditions.

A final remark is that in our formulation,
the compatibility condition and restricted eigenvalue
condition depend on the sparsity $s$ as well as on the active set $S$.
As $S$ is unknown, this means that for a practical
guarantee, the conditions should hold for {\it all} $S$. 
Moreover, one then needs to assume the sparsity $s$ to be
known, or at least a good upper bound needs to be given.
Such strong requirements are the price for practical verifiability.
We however believe that in statistical modeling, non-verifiable
conditions are allowed and in fact common practice. Moreover, our model assumes a
sparse linear ``truth'' with ``true'' active set $S$, only for
simplicity. Without such assumptions, there is no ``true'' $S$, and the
oracle inequality concerns a trade-off between sparse 
approximation and estimation error, see for example \cite{geer08}.


\begin{thebibliography}{24}
\providecommand{\natexlab}[1]{#1}
\providecommand{\url}[1]{\texttt{#1}}
\expandafter\ifx\csname urlstyle\endcsname\relax
  \providecommand{\doi}[1]{doi: #1}\else
  \providecommand{\doi}{doi: \begingroup \urlstyle{rm}\Url}\fi

\bibitem[Bertsimas and Tsitsiklis(1997)]{bertsimas1997introduction}
D.~Bertsimas and J.~Tsitsiklis.
\newblock \emph{{Introduction to linear optimization}}.
\newblock Athena Scientific Belmont, MA, 1997.

\bibitem[Bickel et~al.(2009)Bickel, Ritov, and Tsybakov]{bickel2009sal}
P.~Bickel, Y.~Ritov, and A.~Tsybakov.
\newblock Simultaneous analysis of {L}asso and {D}antzig selector.
\newblock \emph{Annals of Statistics}, 37:\penalty0 1705--1732, 2009.

\bibitem[Bunea et~al.(2007a)Bunea, Tsybakov, and Wegkamp]{Bunea:07a}
F.~Bunea, A.~Tsybakov, and M.~Wegkamp.
\newblock {Aggregation for Gaussian regression}.
\newblock \emph{Annals of Statistics}, 35:\penalty0 1674, 2007a.

\bibitem[Bunea et~al.(2007b)Bunea, Tsybakov, and Wegkamp]{Bunea:07b}
F.~Bunea, A.B. Tsybakov, and M.H. Wegkamp.
\newblock {Sparse Density Estimation with $\ell_1$ Penalties}.
\newblock In \emph{Learning Theory 20th Annual Conference on Learning Theory,
  COLT 2007, San Diego, CA, USA, June 13-15, 2007: Proceedings}, page 530.
  Springer, 2007b.

\bibitem[Bunea et~al.(2007c)Bunea, Tsybakov, and Wegkamp]{Bunea:07c}
F.~Bunea, A.~Tsybakov, and M.~Wegkamp.
\newblock {Sparsity oracle inequalities for the Lasso}.
\newblock \emph{Electronic Journal of Statistics}, 1:\penalty0 169--194, 2007c.

\bibitem[Cai et~al.(2009)Cai, Xu, and Zhang]{cai09}
T.~Cai, G.~Xu, and J.~Zhang.
\newblock On recovery of sparse signals via $\ell_1$ minimization.
\newblock \emph{IEEE Transactions on Information Theory}, 55:\penalty0
  3388--3397, 2009.

\bibitem[Cai et~al.(2009a)Cai, Wang, and Xu]{cai09a}
T.~Cai, L.~Wang, and G.~Xu.
\newblock Shifting inequality and recovery of sparse signals.
\newblock \emph{Preprint}, 2009a.

\bibitem[Cai et~al.(2009b)Cai, Wang, and Xu]{cai09b}
T.~Cai, L.~Wang, and G.~Xu.
\newblock Stable recovery of sparse signals and an oracle inequality.
\newblock \emph{Preprint}, 2009b.

\bibitem[Cand\`es and Plan(2009)]{canpan07}
E.~Cand\`es and Y.~Plan.
\newblock Near-ideal model selection by $\ell_1$ minimization.
\newblock \emph{Annals of Statistics}, 37:\penalty0 2145--2177, 2009.

\bibitem[Cand\`es and Tao(2005)]{candes2005decoding}
E.~Cand\`es and T.~Tao.
\newblock {Decoding by linear programming}.
\newblock \emph{IEEE Transactions on Information Theory}, 51:\penalty0
  4203--4215, 2005.

\bibitem[Cand\`es and Tao(2007)]{candes2007dss}
E.~Cand\`es and T.~Tao.
\newblock {The Dantzig selector: statistical estimation when p is much larger
  than n}.
\newblock \emph{Annals of Statistics}, 35:\penalty0 2313--2351, 2007.

\bibitem[Koltchinskii(2009a)]{koltch09a}
V.~Koltchinskii.
\newblock Sparsity in penalized empirical risk minimization.
\newblock \emph{Annales de l'Institut Henri Poincar\'e, Probabilit\'es et
  Statistiques}, 45:\penalty0 7--57, 2009a.

\bibitem[Koltchinskii(2009b)]{koltch09b}
V.~Koltchinskii.
\newblock The {D}antzig selector and sparsity oracle inequalities.
\newblock \emph{Bernoulli}, 15:\penalty0 799--828, 2009b.

\bibitem[Lounici(2008)]{lounici08}
K.~Lounici.
\newblock {Sup-norm convergence rate and sign concentration property of {L}asso
  and {D}antzig estimators}.
\newblock \emph{Electronic Journal of Statistics}, 2:\penalty0 90--102, 2008.

\bibitem[Meinshausen and B\"uhlmann(2006)]{Meinshausen:06}
N.~Meinshausen and P.~B\"uhlmann.
\newblock {High-dimensional graphs and variable selection with the Lasso}.
\newblock \emph{Annals of Statistics}, 34:\penalty0 1436--1462, 2006.

\bibitem[Meinshausen and Yu(2009)]{meyu09}
N.~Meinshausen and B.~Yu.
\newblock Lasso-type recovery of sparse representations for high-dimensional
  data.
\newblock \emph{Annals of Statistics}, 37:\penalty0 246--270, 2009.

\bibitem[Parter(1961)]{Parter:61}
S.~Parter.
\newblock {Extreme eigenvalues of Toeplitz forms and applications to elliptic
  difference equations}.
\newblock \emph{Transactions of the American Mathematical Society},
  99:\penalty0 153--192, 1961.

\bibitem[van~de Geer(2008)]{geer08}
S.~van~de Geer.
\newblock High-dimensional generalized linear models and the {L}asso.
\newblock \emph{Annals of Statistics}, 36:\penalty0 614--645, 2008.

\bibitem[van~de Geer(2007)]{vandeGeer:07a}
S.~van~de Geer.
\newblock The deterministic {L}asso.
\newblock In \emph{JSM proceedings, {\rm (see also
  http://stat.ethz.ch/research/${\rm research}\_{\rm reports}$/2007/140)}}.
  American Statistical Association, 2007.

\bibitem[Wainwright(2009)]{wainwright09}
M.~Wainwright.
\newblock {Sharp thresholds for high-dimensional and noisy sparsity recovery
  using $\ell_{1}$-constrained quadratic programming ({L}asso)}.
\newblock \emph{IEEE Transactions on Information Theory}, 55:\penalty0
  2183--2202, 2009.

\bibitem[Zhang and Huang(2008)]{zhanghua08}
C.-H. Zhang and J.~Huang.
\newblock The sparsity and bias of the {L}asso selection in high-dimensional
  linear regression.
\newblock \emph{Annals of Statistics}, 36:\penalty0 1567--1594, 2008.

\bibitem[Zhang(2009)]{zhang09}
T.~Zhang.
\newblock Some sharp performance bounds for least squares regression with {L1}
  regularization.
\newblock \emph{Annals of Statistics}, 37:\penalty0 2109--2144, 2009.

\bibitem[Zhao and Yu(2006)]{Zhao:06}
P.~Zhao and B.~Yu.
\newblock {On model selection consistency of Lasso}.
\newblock \emph{Journal of Machine Learning Research}, 7:\penalty0 2541--2563,
  2006.

\bibitem[Zou(2006)]{zou06}
H.~Zou.
\newblock The adaptive {L}asso and its oracle properties.
\newblock \emph{Journal of the American Statistical Association}, 101:\penalty0
  1418--1429, 2006.

\end{thebibliography}

\end{document}